\documentclass[12pt]{amsart}
\usepackage{amsmath,amssymb,amsfonts,latexsym}
\usepackage{mathrsfs}%pour les caracteres en calligraphie
\usepackage[dvips]{geometry}
\usepackage{array}
\DeclareFontFamily{T1}{calligra}{}
\DeclareFontShape{T1}{calligra}{m}{n} {<-> callig15}{}

\newtheorem*{thm}{Theorem}

\newtheorem*{lem}{Lemma}

\newtheorem*{prop}{Proposition}
\theoremstyle{definition}% Pour ne pas avoir les \'enonc\'es en italique.

\newtheorem*{ex}{Example}
\theoremstyle{remark}

\newtheorem*{rem}{Remark}

\numberwithin{equation}{section}
\newcommand{\bT}{\mathbf{Tab}}
\newcommand{\bR}{\mathbf R}
\newcommand{\bS}{\mathbf{S}}
\newcommand{\bB}{\mathbf{B}}
\newcommand{\bG}{\mathbf{G}}
\newcommand{\bM}{\mathbf{M}}
\newcommand{\pr}{^{\prime}}
\newcommand{\prpr}{^{\prime\prime}}
\newcommand{\codim}{{\rm codim\,}}
\newcommand{\Rank}{{\rm Rank\,}}
\newcommand{\sh}{{\rm sh\,}}

\newcommand{\ov}{\overline}
\newcommand{\B}{\mathcal{B}}

\newcommand{\X}{\mathcal{X}}

\newcommand{\N}{\mathcal{N}}

\newcommand{\V}{\mathcal{V}}
\newcommand{\F}{\mathcal{F}}

\newcommand{\Or}{\mathcal{O}}
\newcommand{\Pa}{\mathcal{P}}
\newcommand{\C}{\mathcal{C}}
\newcommand{\D}{\mathcal{D}}
\newcommand{\Bb}{\mathfrak{B}}

\newcommand{\G}{\mathfrak{g}}
\newcommand{\Sl}{\mathfrak{s}\mathfrak{l}}

\newcommand{\gh}{\mathfrak{h}}

\newcommand{\nil}{\mathfrak{n}}

\newcommand{\gN}{\mathfrak{N}}
\newcommand{\Co}{\mathbb{C}}

\newcommand{\Real}{\mathbb R}
\newcommand{\Ze}{\mathbb Z}

\newcommand{\rar}{\rightarrow}

\newcommand{\sr}{\scriptscriptstyle}

\newcommand{\vb}{\vrule height 14pt depth 7pt}
\newcommand{\ts}{\tabskip 4pt}
\newcommand{\vsa}{\noalign{\vskip-7pt}}
\newcommand{\ssa}{\noalign{\vskip -1pt}}

\newcommand{\st}{\subset}
\begin{document}
%\doublespacing%commande pour "double spacing"
%\begin{frontmatter}

\title[B-orbits and patterns]{B-orbits of nilpotent order 2 and link
patterns}%\footnote[1]{(Article soumis)}}%

\author{Anna Melnikov}\thanks{This work was partially supported by
the EEC program RTN-grant "Liegrits"}
\address{Department of Mathematics,
University of Haifa, Haifa 31905, Israel}
\email{melnikov@math.haifa.ac.il}

\begin{abstract}
Let $B_n$ be the group of  upper-triangular  invertible $n\times
n$ matrices and $\X_n$  be the variety of strictly upper
triangular $n\times n$ matrices of nilpotent order 2. $B_n$ acts
on $\X_n$ by conjugation. In this paper we describe geometry of
orbits $\X_n/B_n$ in terms of link patterns.

Further we apply this description to the computations of the
closures of orbital varieties of nilpotent order 2 and intersections
of components of a Springer fiber of nilpotent order 2. In particular we connect our results to the combinatorics of meanders and  Temperley-Lieb algebras.
\end{abstract}
 \maketitle

\section{Introduction}
\subsection{}\label{1.1}
Let $\bM=\bM_n(\Co)$ be an algebra of $n\times n$ matrices over
complex numbers. Let  $\bG=GL_n(\Co)$ be a general linear group.
Consider its action on $\bM$  by conjugation. For $u\in\bM$ let
$\Or_u:=\bG.u:=\{AuA^{-1}\ |\ A\in \bG\}$ denote its orbit.

Let $\gN=\{u\in \bM\ |\ u^k=0\ {\rm for\ some\ } k\}$ be the
nilpotent cone of $\bM.$  Classical Jordan-Gerstenhaber theory
gives a complete combinatorial description of geometry of the
variety of (nilpotent) orbits $\Or_u$ where $u\in \gN$ in terms of
Young diagrams.

Let $\nil=\nil_n\st \bM$ be a subalgebra of strictly
upper-triangular matrices. Let $\bB\st\bG$ be a (Borel) subgroup
of upper-triangular invertible matrices. Consider its action on
$\nil$ by conjugation. For $u\in\nil$ let $\B_u:=\bB.u=\{AuA^{-1}\
|\ A\in\bB\}$ denote its orbit under action of $\bB.$ In general
for $n\geq 6$ the number of such orbits is infinite and there is
no combinatorial theory in the spirit of Jordan-Gerstenhaber
theory describing the geometry of these orbits.

However, if we take a subvariety  of matrices of nilpotent order
2, that is $\X:=\X_n:=\{u\in\nil\ |\ u^2=0\}$ then the number of
$\bB$ orbits in $\X$ is equal to the number of involutions of the
symmetric group $\bS_n.$ We show in the paper that the complete
combinatorial description of geometry of $\bB$ orbits in $\X$ can
be obtained in terms of
 link patterns  in the spirit of
 Jordan-Gerstenhaber theory.

  Further we give two important application
of the technique of link patterns, one in the theory of orbital varieties and another
one in the computations of intersections of the components of a Springer fiber.

\subsection{}\label{1.2}
Let us  first recall in short Jordan-Gerstenhaber theory.

Let $\lambda=(\lambda_1\geq\lambda_2\geq\cdots\geq
\lambda_k>0)\vdash n$ be a partition of $n$ (ordered in decreasing
order). Set $\lambda^*:=\{\lambda^*_1\geq \lambda^*_2\geq \cdots\geq
\lambda^*_l>0\}$ to be the conjugate partition, that is
$\lambda^*_i= \sharp\{j\ |\ \lambda_j\geq i\}.$
\par
 The unique eigenvalue of a nilpotent matrix is 0, so that Jordan
form of $u\in \gN$ is defined completely by the lengths of its
Jordan blocks, which is defined in turn by some partition of $n$.
Since each nilpotent orbit has a unique Jordan form up to the
order of Jordan blocks we have a bijection between nilpotent
orbits and ordered partitions.  We can write $J(u):=\lambda$ and
$\Or_\lambda:=\Or_u$ where $\lambda$ is the corresponding
partition.

 Let $\ov V$ denote the closure of variety $V$ with
respect to Zariski topology. The geometry of $\ov\Or_\lambda$ is
described combinatorially in terms of partitions. We begin with
the formula of the dimension of an orbit:
$$\dim \Or_\lambda=n^2-\sum\limits_{i=1}^l(\lambda_i^*)^2.\eqno{(*)}$$

Our next goal is  a combinatorial description of the closure of an
orbit.

 A partial order on partitions (called a dominance
order) is defined as follows. Let $\lambda,\mu\vdash n$ be ordered
in decreasing order $\lambda=(\lambda_1\geq \lambda_2\geq \ldots
\lambda_k>0)$ and $\mu=(\mu_1\geq \mu_2\geq\ldots\mu_l>0)$. Set
$\lambda\geq \mu$ if for any $j\ :\ 1\leq j\leq \min(k,l)$ one has
$\sum\limits_{i=1}^j\lambda_i\geq \sum\limits_{i=1}^j\mu_i.$

Let $\coprod$ denote a disjoint union. By a theorem of M.
Gerstehaber (cf. \cite{He} for example)
$$\overline\Or_\lambda=\coprod\limits_{\mu\leq \lambda}\Or_\mu.$$

Given a partition $\lambda=(\lambda_1\geq \ldots\geq \lambda_k>0)$
we define the corresponding Young diagram $D_\lambda$ of $\lambda$
to be an array of $k$  rows of cells starting on the left with the
$i$-th row containing $\lambda_i$ cells.
\par
Young diagrams are a very convenient visualization of the partitions
with respect to nilpotent orbits. Indeed, the dual partition used in
$(*)$ is simply the list of lengths of the columns of the
corresponding Young diagrams. Also, $\mu<\lambda$ if $D_\mu$ is
obtained from $D_\lambda$ by pushing down some cells. In particular,
the cover of a given partition (that is all the maximal partitions
strictly smaller than the given one) with respect to the dominance
order  is described easily in the terms of Young diagrams.

Indeed, let $\mu=(\mu_1,\ldots)$ be in the cover of
$\lambda=(\lambda_1,\ldots).$ Then $\mu$ is obtained from $\lambda$
in one of two ways:
\begin{itemize}
\item[(i)] There exists $i$ such that $\lambda_i-\lambda_{i+1}\geq 2$. \ Then
$\mu_j = \lambda_j$ for $j\neq i, \ i+1$  and
$$ \mu_i = \lambda_i-1, \ \ \mu_{i+1} = \lambda_{i+1}+1 \ . $$
In this case $\dim\Or_\mu=\dim\Or_\lambda-2.$
\medskip
\item[(ii)] There exists $i$ such that
            $ \lambda_{i+1}=\lambda_{i+2}=\cdots=\lambda_{i+k}= \lambda_i-1$
for  some  $k\geq 1$  and $\lambda_{i+k+1} = \lambda_i-2.$ Then
$\mu_j = \lambda_j$ for $j\neq i, \ i+k+1$ \ and
$$ \mu_i = \lambda_i-1, \ \ \mu_{i+k+1} = \lambda_i-1 \ . $$
In this case $\dim\Or_\mu=\dim\Or_\lambda-2(k+1).$
\end{itemize}
In particular the cover of a nilpotent orbit  is not equidimensional
in general.

 The above result can be
described by Young diagrams as follows.
\par
In the first case $D_\mu$ is obtained from $D_\lambda$ by pushing
one box down  one row (and possible across several columns). For
example,
$$ D_\lambda =
\vcenter{ \halign{& \hfill#\hfill \tabskip4pt\cr
\multispan{11}{\hrulefill}\cr
 \ssa \vb & \quad & \vb & \quad & \vb &  \quad & \vb & \quad & \vb & \quad &\ts\vb\cr
 \vsa \multispan{11}{\hrulefill}\cr
 \ssa \vb & \ \ & \vb & \ \ &\vb & \ \  & \vb  & \ \  & \vb &X&\ts\vb\cr
  \vsa
\multispan{11}{\hrulefill}\cr \ssa \vb & \ & \ts\vb \cr \vsa
\multispan{3}{\hrulefill}\cr \ssa \vb  & \ &\ts\vb\cr \vsa
\multispan{3}{\hrulefill}\cr}}\ , \qquad
 D_\mu =
\vcenter{ \halign{& \hfill#\hfill \tabskip4pt\cr
\multispan{11}{\hrulefill}\cr
 \ssa \vb & \quad & \vb &  \quad & \vb &  \quad & \vb
& \quad & \vb & \quad &\ts\vb\cr
 \vsa \multispan{11}{\hrulefill}\cr
\ssa \vb & \ & \vb & \ & \vb &  \ &  \vb & \ &\ts\vb\cr \vsa
\multispan{9}{\hrulefill}\cr \ssa \vb & \ & \vb &X&\ts\vb\cr \vsa
\multispan{5}{\hrulefill}\cr \ssa \vb  & \ &\ts\vb\cr \vsa
\multispan{3}{\hrulefill}\cr}}\ . $$

\par
In the second case diagram $D_\mu$ is obtained from $D_\lambda$ by
pushing one box across one column (and possible down several
rows). For example,
$$ D_\lambda =
\vcenter{ \halign{& \hfill#\hfill \tabskip4pt\cr
\multispan{9}{\hrulefill}\cr
 \ssa \vb & \quad & \vb &  \quad & \vb&\quad & \vb &X &\ts\vb\cr
 \vsa \multispan{9}{\hrulefill}\cr
 \ssa\vb & \ & \vb & \ & \vb  & \ &\ts\vb\cr
 \vsa\multispan{7}{\hrulefill}\cr
 \ssa\vb & \ & \vb & \ & \vb  & \ &\ts\vb \cr
 \vsa\multispan{7}{\hrulefill}\cr
 \ssa\vb & \ & \vb & \ & \ts\vb \cr
 \vsa\multispan{5}{\hrulefill}\cr
 \ssa \vb  & \ &\ts\vb\cr
 \vsa\multispan{3}{\hrulefill}\cr}}\ , \qquad
 D_\mu =
\vcenter{ \halign{& \hfill#\hfill \tabskip4pt\cr
 \multispan{7}{\hrulefill}\cr
 \ssa \vb & \quad & \vb &  \quad & \vb & \quad & \ts\vb\cr
 \vsa \multispan{7}{\hrulefill}\cr
 \ssa \vb & \ &\vb &  \ &  \vb & \ & \ts\vb\cr
 \vsa \multispan{7}{\hrulefill}\cr
 \ssa \vb & \ &\vb &  \ &  \vb & \ & \ts\vb\cr
 \vsa \multispan{7}{\hrulefill}\cr
\ssa \vb & \ & \vb &  \ &  \vb &X&  \ts\vb\cr \vsa
\multispan{7}{\hrulefill}\cr \ssa \vb  & \ &\ts\vb\cr \vsa
\multispan{3}{\hrulefill}\cr}}\ .$$

\subsection{}\label{1.3}
 In the paper we show that link patterns
play the same role for the combinatorial description of geometry of
$\bB$ orbits in $\X$ as Young diagrams for the description of
geometry of nilpotent orbits. First we describe the bijection
between $\bB$ orbits in $\X$ and link patterns.

Let $\bS_n^2:=\{\sigma\in\bS_n\ |\ \sigma^2=Id\}$  be the subset
of involutions of the symmetric group $\bS_n$. For any $\sigma\in
\bS_n^2$ let $N_\sigma$ be the matrix obtained from the
representation matrix of $\sigma$ by erasing the lower-triangular
part (including the main diagonal), that is
$$(N_\sigma)_{i,j}:=\left\{\begin{array}{ll} 1&{\rm if}\ i<j\
{\rm and}\
\sigma(i)=j;\\
0&{\rm otherwise.}\\
\end{array}\right.$$
One can see at once that $N_\sigma\in\X$ for any $\sigma\in\bS_n^2.$
Moreover by \cite[2.2]{Mx2} one has
$$\X=\coprod\limits_{\sigma\in\bS_n^2}\B_{N_{\sigma}}.$$
Thus, $\bB$ orbits in $\X$ are labelled by involutions of $\bS_n.$
Put $\B_\sigma:=\B_{N_\sigma}.$

Given $\sigma\in\bS_n^2$ the corresponding link pattern $P_\sigma$
is an array of $n$ points on a (horizontal) line where  points
$i\ne j$ are connected by an arc if $\sigma(i)=j.$ We will call
such points end points of an arc. A point $i$ satisfying
$\sigma(i)=i$ is called a fixed point of $P_\sigma$. For example,
for $\sigma=(1,3)(2,6)(4,7)\in \bS_7^2$ one has
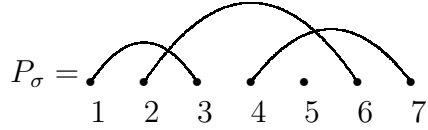
\begin{center}
\begin{picture}(100,80)
\put(-20,40){$P_\sigma=$}
\multiput(10,40)(20,0){7}%
{\circle*{3}}
 \put(10,25){1}
 \put(30,25){2}
 \put(50,25){3}
 \put(70,25){4}
 \put(90,25){5}
 \put(110,25){6}
 \put(130,25){7}
 \qbezier(10,40)(30,70)(50,40)
 \qbezier(30,40)(70,100)(110,40)
 \qbezier(70,40)(100,80)(130,40)
\end{picture}
\end{center}
The only fixed point of $P_\sigma$ is $5.$

\begin{rem}I found the notion of link pattern in the papers in
combinatorics and mathematical physics \cite{Z}, \cite{Kn},
\cite{Z-J}. In general link patterns considered there are patterns
with maximal possible number of  arcs and a link pattern is a
circle with n points rather than a line.  Respectively the authors
use a word "arch" not "arc". I adopted the notion of a link
pattern from them, however it seems to me that the word "arc" is
more appropriate in the representation given here.
\end{rem}
\subsection{}\label{1.4}
Let us explain the role of link patterns in the description of
geometry of $\bB$ orbits in $\X$.

Let $\B_\sigma$ be an orbit in $\X$ and $P_\sigma$ be the
corresponding link pattern. We begin with formula of
$\dim\B_\sigma$ in terms of a link pattern. Given
$\sigma=(i_1,j_1)\ldots(i_k,j_k)\in \bS_n^2$ (in cyclic form).
 \begin{itemize}
 \item{} Put $\ell(P_\sigma):=k$ to be the number of arcs
in $P_\sigma$ (we will call it the length of a pattern);
\item{} Arcs $(i,j),(i',j')$ where $i<i'$  of $P_\sigma$ are called intersecting if
$i<i'<j$ and $j'>j.$ We say that arcs $(i,j), (i',j')$ have a
crossing in this case. Put $c(P_\sigma)$ to be the number of
crossings of arcs in $P_\sigma;$
\item{} Let $\{p_s\}_{s=1}^{n-2k}$ be the set of fixed points
of $P_\sigma$. For $s:\ 1\leq s\leq n-2k$ put $f_{p_s}(P_\sigma)$ to
be the number of arcs $(i_t,j_t)$ over $p_s$ (that is such that
$i_t<p_s<j_t$) and $f(P_\sigma):=\sum_{s=1}^{n-2k}
f_{p_s}(P_\sigma).$
\end{itemize}
 In the example above
$\ell(P_\sigma)=3$, $c(P_{\sigma})=2$, and
$f(P_\sigma)=f_{5}(P_\sigma)=2.$

As we show in \ref{3.1} one has
$$\dim
\B_\sigma=\ell(P_\sigma)\cdot(n-\ell(P_\sigma))-c(P_\sigma)-f(P_\sigma).$$
Thus, in our example $\dim \B_\sigma=3\cdot 4-2-2=8.$

Note that $\bS_n^2$ can be stratified by the length of link
patterns. Put $\bS_n^2(k):=\{\sigma\in\bS_n^2\ |\
\ell(P_\sigma)=k\}.$ One has $\bS_n^2=\coprod\limits_{k\leq{n\over
2}}\bS_n^2(k).$

 Note
that this stratification partitions $\bB$ orbits in $\X$ into sets
belonging to the same nilpotent orbit. Indeed, put
$\Or_\sigma:=\Or_{N_\sigma},$ then $\Or_\sigma=\Or_{\sigma'}$  iff
$\ell(P_\sigma)=\ell(P_{\sigma'}).$ This is because a nilpotent
orbit of nilpotent order 2 is completely defined by the rank of
its representing matrix and $\Rank(N_\sigma)=\ell(P_\sigma).$

\subsection{}\label{1.5}

Further we describe $\ov\B_\sigma$  in terms of $P_\sigma$.

 For $i,j\ :\ 1\leq i\leq j\leq
 n$  let $[i,j]:=\{s\}_{s=i}^{j}$ be the set of integer points of a segment $[i,j]$. For an arc $(i_s,j_s)\in P_\sigma$ we
 put
 $(i_s,j_s)\in [i,j]$
 if $i\leq i_s,j_s\leq j.$
  For $i,j\ :\ 1\leq i<j\leq n$ let $(R_\sigma)_{i,j}$
be the number of  arcs in $P_\sigma$ belonging to $[i,j]$. Just
for convenience we can represent $R_\sigma$ as an $n\times n$
matrix with integer non negative entries by setting
$(R_\sigma)_{i,j}=0$ for any $i\geq j.$ In our example

$$R_{\sigma}=\left(\begin{array}{ccccccc}
0&0&1&1&1&2&3\\
0&0&0&0&0&1&2\\
0&0&0&0&0&0&1\\
0&0&0&0&0&0&1\\
0&0&0&0&0&0&0\\
0&0&0&0&0&0&0\\
0&0&0&0&0&0&0\\
\end{array}\right)$$

Define a partial order on $\bM_n(\Real)$ by $A\succeq B$ if for any
$i,j\ :\ 1\leq i,j\leq n$ one has $(A)_{i,j}\geq (B)_{i,j}.$ This
order induces a partial order on involutions by putting
$\sigma\succeq \sigma'$ if $R_\sigma \succeq R_{\sigma'}.$ Then, by
\ref{2.4} and \ref{3.2} one has
$$\overline\B_\sigma=\coprod\limits_{\sigma'\preceq\sigma}\B_{\sigma'}.$$

\subsection{}\label{1.6}
Let us apply the previous result to the cover of a given $\sigma$
with respect to order $\preceq$. We describe the cover  in terms of
link patterns only.

First, notice that nilpotent orbits of nilpotent order 2 are
ordered linearly with respect to a dominance order. Indeed, given
$\lambda, \mu$ corresponding to nilpotent orbits of nilpotent
order 2  then $\lambda^*=(n-k,k)$ and $\mu^*=(n-m,m)$. One has
$\mu\leq\lambda$ iff $m\leq k$ otherwise $\mu>\lambda.$ Thus, for
any nilpotent orbits $\Or,\Or'$ of nilpotent order 2 one has
either $\ov\Or'\subset \ov\Or$ or $\ov\Or'\supset\ov\Or.$ In
particular the cover of partition $(n-k,k)^*$ consists  of the
unique partition $(n-k+1,k-1)^*$ ( it is of type (ii) described in
\ref{1.2}).

Let $\C(\sigma)$ denote the cover for a given $\sigma\in \bS_n^2$
with respect to order $\preceq.$ In \ref{3.9} we show that
$\sigma'\in \C(\sigma)$ iff $\codim_{\ov\B_\sigma}\B_{\sigma'}=1$ so
that the cover with respect to order $\preceq$ is equidimensional.

 However, from geometric point of view there are two
interesting subsets of elements smaller than $\sigma.$ The first one
is $\N(\sigma)$  -- the subset of maximal elements smaller than
$\sigma$ in the set of patterns of smaller length. The second one is
$\D(\sigma):=\C(\sigma)\cap\bS_n^2(\ell(P_\sigma)).$

Let $\sigma\in\bS_n^2(k).$ Geometric meaning of $\N(\sigma)$ is
provided by
$$\ov\B_\sigma\bigcap(\ov\Or_\sigma\setminus
\Or_\sigma)=\ov\B_\sigma\bigcap\ov\Or_{(n-k+1,k-1)^*}=
\bigcup\limits_{\sigma'\in \N(\sigma)}\ov \B_{\sigma'}.$$

As for  $\D(\sigma)$ it helps  to understand the structure of
$\ov\B_\sigma$ inside of $\Or_{\sigma}$.
 In particular for $\sigma,\sigma'\in\bS_n^2(k)$ such that
 $\dim \B_\sigma=\dim \B_{\sigma'}$
one has $\codim_{\ov\B_\sigma}(\ov\B_\sigma\cap\ov\B_{\sigma'})=1$
if and only if $\D(\sigma)\cap \D(\sigma')\ne\emptyset.$

 We give the combinatorial description in terms of link
patterns of both $\N(\sigma)$ and $\D(\sigma)$. The description of
$\N(\sigma)$  in terms of link patterns is very simple. We call an
arc of $P_\sigma$ {\bf external} if there is no arc over it, that
is $(i,j)\in P_\sigma$ is external if for any arc
$(i',j')\ne(i,j)$ in $P_\sigma$ one has $(i,j)\not\in[i',j'].$

All the elements of $\N(\sigma)$ are obtained from $\sigma$ by
deleting an external arc. So the number of elements in $\N(\sigma)$
is equal to the number of external arcs in $\sigma$ and all the
elements in $\N(\sigma)$ are in $\bS_n^2(\ell(P_\sigma)-1).$

Note that in general
\begin{itemize}
\item{} the set $\{ \B_{\sigma'}\ |\ \sigma'\in
N(\sigma)\}$ is not equidimensional;
 \item{} $\N(\sigma)$ is {\bf not} a subset of $C(\sigma).$
\end{itemize}
$\D(\sigma)$ consists of 3 possible types of elements. We describe
them in \ref{3.4}-\ref{3.6}.

\subsection{}\label{1.6a}
 All the results mentioned above are based on \cite{Mo2}
and we use heavily the results of this paper in our proofs. The
advantage of the representation of involutions by link patterns is
that  the complex combinatorics of \cite{Mo2} becomes easy and
natural in this language.

 The fact that $\codim_{\ov\B_\sigma}\B_{\sigma'}=1$ for any
$\sigma'\in \C(\sigma)$ was shown already in \cite{M-PII}.
However, the proof provided there involves a lot of computations
and here using the technique of link patterns we prove it in a few
lines.

\subsection{}\label{1.7} Let us use the technique of link patterns
 to study  orbital varieties of nilpotent
order 2 and  components of Springer fibers of nilpotent order 2. To
do this we consider  link patterns without intersecting arcs and
without fixed points under the arcs. They are the most simple link
patterns from the combinatorial point of view on one hand and the
most important objects for the applications to representation theory
on the other hand.

By \ref{1.4} $\dim\B_\sigma=(n-k)k$ (that is $\B_\sigma$ is of
maximal possible dimension) iff $P_\sigma$ is such link pattern.
Note that by \ref{1.2} $\dim
\Or_{(n-k,k)^*}=n^2-k^2-(n-k)^2=2k(n-k).$ Since in general
$\dim\B_u\leq 0.5\dim\Or_u$ for $u\in\nil$ we get that $\B_\sigma$
corresponding to link patterns without intersecting arcs and without
fixed points under the arcs are $\bB$ orbits of maximal possible
dimension. Moreover, their closures in  $\Or_\sigma$ are Lagrangian
subvarieties of $\Or_\sigma$, called orbital varieties.

Let us give a brief general description of an orbital variety (cf.
\ref{2.2} for the details). Let $\Or$ be a nilpotent orbit in a
semisimple Lie algebra $\G$. Let us fix some triangular
decomposition $\G=\nil\oplus\gh\oplus \nil^-$. The irreducible
components of $\Or\cap\nil$ are called orbital varieties
associated to $\Or.$ Note that $\Or$ is a symplectic manifold and
as it was shown by N. Spaltenstein, R. Steinberg and A. Joseph an
orbital variety is a Lagrangian subvariety of $\Or.$ According to
orbit method Lagrangian subvarieties should play an important role
in the representation theory.

Let us explain the role of orbital varieties in $\G=\Sl_n$
considered in the paper. In this case orbital varieties associated
to $\Or_\lambda$ are parameterized by standard Young tableaux
corresponding to Young diagram $D_\lambda.$ Let $\bT_\lambda$ be
the set of standard Young tableaux corresponding to $D_\lambda.$
For $T\in\bT_\lambda$ we denote by $\V_T$ the corresponding
orbital variety associated to $\Or_\lambda.$

By a result of A. Joseph primitive ideals of the enveloping algebra
$U(\G)$ corresponding to a highest weight module with an integral
weight are in bijection with standard Young tableaux. Let  $I_T$ be
the primitive ideal corresponding to $T\in\bT_\lambda.$ As it is
shown in \cite{M} an orbital variety closure $\overline \V_T$ is the
associated variety of $I_T.$  As those they play a key role in the
theory of primitive ideals. The details on the theory above can be
found in \cite{J}.

Let $S,T$ be Young tableaux with two columns. By \cite{Mrang2}
$I_T\subset I_S$ if and only if $\ov \V_T\supset \ov \V_S.$
Further in \cite{Mo2} it was shown that for $T\in \bT_{(n-k,k)^*}$
one has that the closure of $\V_T$ is a union of orbital varieties
and the combinatorial description of this closure in terms of
Young tableaux was obtained. However the combinatorial description
in \cite[4.3]{Mo2} is too complex to be satisfactory. On the other hand,
as it was noted already in \cite{Msmith} each orbital variety $\V_T$ of nilpotent order two
admits the unique $\bB$ orbit. Let us denote the corresponding involution by $\sigma_T$. In \ref{4.2}
we simplify this description using link patterns and get a simple
answer in terms of $\N(\sigma_T)$.

\subsection{}\label{1.8}
 As it was shown in \cite{M-P} the bijection
between orbital varieties associated to a nilpotent orbit
$\Or_\lambda$ and components of Springer fiber $\F_x$ where  $x\in
\Or_\lambda$ is extended to the intersections: Let $S,T\in
\bT_\lambda$ and let $\V_T,\V_S$ be orbital varieties associated to
$\Or_\lambda$ and $\F_T,\F_S$ be the corresponding components of
$\F_x$. The number of irreducible components of $\V_T\cap\V_S$ and
their codimensions are equal to the number of irreducible components
of $\F_T\cap\F_S$ and their codimensions. Via this bijection orbital
varieties play a key role in Springer theory.

 This bijection together with the theory of $\bB$ orbits of nilpotent
order 2 gave the opportunity to compute for the first time the
intersections of the components of Springer fiber when these
components are singular. These intersections were considered in
\cite{M-PI}, \cite{M-PII}. Here we use link patterns to simplify
and extend these results. We start with a general theory of
intersections of the closures of $\bB$ orbits of nilpotent order 2 belonging to
$\Or_{(n-k,k)^*}.$
 Let $\sigma,\sigma'\in\bS_n^2(k)$.
 The general picture of
$\ov\B_\sigma\cap\ov\B_{\sigma'}\cap\Or_{(n-k,k)^*}$
  was considered in \cite{M-PI},
\cite{M-PII}. As it is shown there these intersections can be very
complex. In particular,there  are examples in these papers of a reducible intersection  of codimension 1 and of a reducible intersection of higher codimension which is
not equidimensional ( in \ref{4.3} we give an example of an intersection of codimension 1 which is not equidimensional). To compute
$\ov\B_\sigma\cap\ov\B_{\sigma'}\cap\Or_{(n-k,k)^*}$  we use matrix
$R_{\sigma,\sigma'}$  defined by
$$(R_{\sigma,\sigma'})_{i,j}:=\min\{(R_{\sigma})_{i,j},\ (R_{\sigma'})_{i,j}\}$$
As it was shown in \cite[5.6, 5.7]{M-PI} this intersection is non
empty and it is a union of all $\B_{\sigma^{\prime\prime}}$ where
$\sigma^{\prime\prime}\in \bS_n^2(k)$  are such that
$R_{\sigma^{\prime\prime}}\preceq R_{\sigma,\sigma'}$. In
particular this intersection  is irreducible iff there exists
$\sigma^{\prime\prime}\in\bS_n^2(k)$ such that
$R_{\sigma,\sigma'}=R_{\sigma^{\prime\prime}}.$

To simplify the results on intersections we pass from link patterns
to so called (generalized) meanders. To do this we draw two link
patterns on the same line of points in such a way that the arcs of
the first link pattern are drown  upward and  the arcs of the second
one are drown downward.

 For example, consider $\sigma=(1,3)(2,5)(4,7),\ \sigma'=(1,4)(3,5)(6,7)\in \bS_7^2.$
Their generalized meander $M_{\sigma,\sigma'}$  is
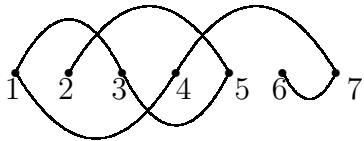
\begin{center}
\begin{picture}(100,80)
\multiput(10,40)(20,0){7}%
{\circle*{3}}
 \put(6,30){1}
 \put(26,30){2}
 \put(46,30){3}
 \put(70,30){4}
 \put(92,30){5}
 \put(106,30){6}
 \put(134,30){7}
 \qbezier(10,40)(30,80)(50,40)
 \qbezier(30,40)(60,90)(90,40)
 \qbezier(70,40)(100,90)(130,40)

\qbezier(10,40)(40,-10)(70,40)
 \qbezier(50,40)(70,0)(90,40)
 \qbezier(110,40)(120,20)(130,40)
\end{picture}
\end{center}
These generalized meanders are useful for computing
$R_{\sigma,\sigma'}.$

Taking a meander built of 2  link patterns corresponding to
maximal $\bB$ orbits in $\Or_{(n-k,k)^*}$  we get a classical
meander (that is without crossing arcs) without fixed points under
any (upward or downward) arc. As it is shown in \cite{M-PII}  for
$S,T\in \bT_{(n-k,k)^*} $ if  $\codim_{\V_T}(\V_S\cap\V_T)=1$ then
$\V_S\cap\V_T$ is irreducible.  The combinatorial description of
pairs $S,T$ with intersection of codimension 1 is provided there.

However, the consideration of these intersections in terms of
meanders first of all simplifies the description and the proofs and
further connects our results to Temperley-Lieb representations of
Hecke algebras and Kazhdan-Lusztig data.

Let us formulate the result. A meander $M_{\sigma_T,\sigma_S}$ is
called even if any its connected part consists of even number of
arcs. Otherwise $M_{\sigma_T,\sigma_S}$ is called  odd. The connected part consisting of even number
of arcs may be either a closed path (we call it a loop) or an open
path (we call it an even interval). As we show in \ref{4.5} for $S,T\in
\bT_{(n-k,k)^*} $ one has $\codim_{\V_T}(\V_S\cap\V_T)=1$ if and
only is $M_{\sigma_S,\sigma_T}$ is an even meander with $k-1$ loops.
The irreducibility of such intersections is obvious from the combinatorics of
the meanders.

B.W. Westbury in \cite{W} computed the inner product matrix for the
Kazhdan-Lusztig basis for a two row shape tableaux in terms of
Temperley-Lieb cup diagrams. Note that for $T\in\bT_{(n-k,k)}$ the
corresponding cup diagram defined in \cite{W} is exactly
$P_{\sigma_{T^t}}$ where $T^t\in\bT_{(n-k,k)^*}$ is a transposed
tableau. Given $S,T\in\bT_{(n-k,k)}$ let $c_T,c_S$ be the
corresponding Kazhdan-Lusztig basis vectors from \cite{W}. As it is
shown in
 \cite[7]{W} their inner
product $< c_S,c_T>$ is always either 0 or $(t+t^{-1})^r$ where
$0\leq r\leq k-1.$ Moreover, restating his results in the language
of meanders we get $< c_S,c_T>=(t+t^{-1})^{r}$ if and only if the
corresponding meander $M_{\sigma_S,\sigma_T}$ is  an even meander with $r$
loops.

 F.Y.C. Fung in \cite{F} computed the
 intersections of the components of a Springer fiber corresponding to
 a nilpotent element with two Jordan blocks, that is $\F_T\cap \F_S$
 where $T,S\in\bT_{(n-k,k)}$ . In this case
 $\F_T$ is
 nonsingular and as he showed all the intersections are either
 irreducible or empty. Using the technique similar to those used in \cite{W}
 he showed that $\F_S\cap\F_T=\emptyset$
 if and only if $< c_T,c_S>=0$ and if $< c_T,c_S>\ne 0$ then
 $\codim_{\F_T}(\F_T\cap \F_S)=k- \deg< c_T,c_S>.$

Consider $S,T\in\bT_{(n-k,k)}.$
 Combining the above we get that $\codim_{\F_T}(\F_T\cap \F_S)=1$ if
 and only if $\codim_{\F_{T^t}}(\F_{T^t}\cap \F_{S^t})=1$ and this
 happens iff $<c_S,c_T>=(t+t^{-1})^{k-1}.$

However the pictures of the intersections of higher codimension in two-row and two-column cases are very different.
 Note that for $\sigma,\sigma'\in\bS_n^2(k)$ one has
 $\ov\B_\sigma\cap\ov\B_{\sigma'}\cap\Or_{(n-k,k)^*}\ne\emptyset$, thus, in particular
 for an odd meander $M_{\sigma_S,\sigma_T}$ one has $<c_{S^t},c_{T^t}>=0$
 and $\F_T\cap\F_S\ne\emptyset.$ On the other hand, for even meanders with less than $k-1$ loops it can occur that
  $\codim_{\F_T}(\F_T\cap \F_S)> \codim_{\F_{T^t}}(\F_{T^t}\cap \F_{S^t})$ as we show
  by example in \ref{4.7}. The general theory of intersections of higher codimension
  is too complex from combinatorial point of view and not that interesting from geometric point of view so we do not discuss it here.

  Also, the intersections of higher codimension are almost always reducible. Again, the general theory is too involved combinatorially to explain it here so we only provide in \ref{4.8} a very simple combinatorial sufficient condition for the reducibility of an intersection which shows that in general the intersections of higher codimension
  are reducible.

\vskip 0.2 cm
 In the end of the paper one can find the index of notation in which
symbols appearing frequently are given with the subsection where
they are defined. I hope that this will help the reader to find
his way through the paper.

 \vskip 0.2 cm
 {\bf Acknowledgements.} I am very grateful to A.
 Elashvili for the discussions on biparabolic algebras and the role
 of meanders in their description. He was first to introduce the
 world of meanders to me. Through these discussions I realized the
 role that they (and link patterns) play in the description of $\bB$
 orbits of nilpotent order 2 and intersections of the components of
 a Springer fiber of nilpotent order 2. I am also very grateful to A. Joseph
 for many fruitful discussions during various stages of this work.

\section{$\bB$ orbits in $\X$ and involutions}

\subsection{}\label{2.1}
Put $\Bb_{(n-k,k)}=\coprod\limits_{\sigma\in\bS_n^2(k)}\B_\sigma.$
As we noted in \ref{1.4} $\Bb_{n-k,k}=\X\cap\Or_{(n-k,k)^*}$ so that
$\X=\coprod\limits_{k\leq {n\over 2}}\Bb_{(n-k,k)}.$

We start with the computation of the dimension of $\B_\sigma.$
  All over this paper we write a cycle in round brackets
$(i,j)$ when the entries  are ordered in the increasing order, i.e.
$i<j$. We write it in double round brackets $((i,j))$ when the entries $i\ne j$ are
not ordered.
 Write $\sigma\in \bS_n^2(k)$ as a product of disjoint cycles
of length 2, that is $\sigma=(i_1,j_1)(i_2,j_2)\ldots(i_k,j_k)$
where $i_s<j_s$ and $\{i_s,j_s\}\cap\{i_t,j_t\}=\emptyset$ for any
$1\leq s\ne t\leq k.$ We write $(i,j)\in \sigma$ or $(i,j)\in
P_\sigma$ if $(i,j)=(i_s,j_s)$ for some $s\ :\ 1\leq s\leq t.$ For
$(i,j)\in\sigma$ put

$$q_{(i,j)}(\sigma):=\#\{i_p<i\ | \ j_p<j\}+
\#\{j_p\ |\ j_p<i\}. \eqno{(*)}
$$

\begin{ex} Take $\sigma=(1,6)(3,4)(5,7)\in\bS_7^2(3).$ Then $\ell(P_\sigma)=3$ and
$q_{(1,6)}(\sigma)=0,\ q_{(3,4)}(\sigma)=0,\
q_{(5,7)}(\sigma)=2+1=3.$
\end{ex}
\par
By \cite[3.1]{Mx2} one has
\begin{thm}\label{dim}
 For $\sigma=(i_1,j_1)(i_2,j_2)\ldots(i_k,j_k)\in \bS_n^2(k)$
   one has
$$\dim \B_\sigma=kn-\sum\limits_{s=1}^k(j_s-i_s)-\sum\limits_{s=1}^k q_{(i_s,j_s)}(\sigma).$$
\end{thm}
\subsection{}\label{2.2}
What is the maximal  dimension of an orbit in $\Bb_{(n-k,k)}$? We
would like to provide a detailed answer to this question. To do
this we consider in detail orbital varieties defined briefly in
\ref{1.7}.

Let $\bG$ be a connected semisimple finite dimensional complex
algebraic group and let $\G={\rm Lie}(\bG)$ be the corresponding
semisimple Lie algebra. Let $\G=\nil\oplus\gh\oplus\nil^-$ be its
triangular decomposition.  For $u\in\nil$ let $\Or_u=\bG.u$ be its
adjoint (nilpotent) orbit. Note that by the identification of $\G$
with $\G^*$ through the Killing form $\Or_u$ is provided a
symplectic structure.

Consider $\Or_u\cap\nil.$ This variety is reducible and its
irreducible components are called orbital varieties associated to
$\Or_u$. Let
 $\V$ be an orbital variety associated to
$\Or_u.$  By ~\cite{Sp2} and ~\cite{St1} one has $\dim \V={1\over
2}\dim\Or_u.$ ( Moreover as it was pointed out in ~\cite{J} this
implies that an orbital variety is a Lagrangian subvariety of the
nilpotent orbit it is associated to.)

In particular, if $\bG=SL_n$ nilpotent orbits are described by
Young diagrams as we explained in \ref{1.2}. In this case there
exists a very nice combinatorial characterization of orbital
varieties in terms of Young tableaux. Recall that a Young tableau
$T$ associated to  Young diagram $D_\lambda$ is obtained by
filling the boxes of $D_\lambda$ with numbers $1,\ldots,n$ so that
the numbers increase in rows from left to right and in columns
from top to bottom. Given a Young tableau $T$ associated to
$D_\lambda$ its shape is defined to be $\lambda$ and denoted by
$\sh(T).$ Given $u\in \nil\cap \Or_\lambda$ its Young diagram is
again defined to be $\lambda$ and denoted by $D_n(u)$, or simply
$D(u).$ Now consider canonical projections $\pi_{1,n-i}:\nil_n\rar
\nil_{n-i}$ acting on a matrix by deleting the last $i$ columns
and the last $i$ rows. Set $D_{n-1}(u):=D(\pi_{1,n-1}(u)),
D_{n-2}(u):=D(\pi_{1,n-2}(u)),\ldots, D_1(u):=D(\pi_{1,1}(u))$ and
$\phi(u):=(D_1(u),D_2(u),\ldots,D_n(u)).$ Put $1$ into the unique
box of $D_1(u).$ For any $i\ : 1\leq i\leq n-1$ note that
$D_{i+1}(u)$ differs from $D_i(u)$ by a single "corner" box. Put
$i+1$ into this box. This gives a bijection from the set of the
chains $(D_1(u),D_2(u),\ldots,D_n(u))$ to the set of standard
Young tableaux $T$ of shape $D_n(u).$ In other words, we view a
standard Young tableau as a chain of Young diagrams. For
$u\in\nil$ put $\varphi(u):=T$ if $T$ corresponds to $\phi(u)$
under this bijection. Set $\nu_{\sr T}:=\{u\in \nil\ |\
\varphi(u)=T\}.$
\par
By Spaltenstein ~\cite{Sp} orbital varieties associated to
$\Or_\lambda$ are parameterized by standard Young tableaux of shape
$\lambda$ as follows. Let $\{T_i\}$ be the set of Young tableaux of
shape $\lambda.$ Set $\V_{T_i}:=\ov\nu_{\sr T_i}\cap \Or_\lambda.$
Then  $\{\V_{T_i}\}$ is the set of orbital varieties associated to
$\Or_\lambda.$

Note that in particular by this construction $\dim \B_u\leq
{1\over 2} \dim\Or_u.$ In general this inequality is strict.
Moreover, obviously if $\dim \B_u={1\over 2} \dim\Or_u$ then
$\B_u$ is dense in the corresponding orbital variety.
Unfortunately, for $n\geq 6$ the vast majority of orbital
varieties do not admit a dense $\bB$ orbit. However, orbital
varieties of nilpotent order two always admit a unique dense
orbit. That gives the answer to the question posed above.

To formulate this answer let us consider the corresponding Young
tableaux in detail. Let $\bT_{(n-k,k)^*}$ be the set of Young
tableaux of shape $(n-k,k)^*.$
 For $T\in \bT_{(n-k,k)^*}$ set $T=(T_1,T_2),$ where
$T_1=\left(\begin{array}{c}a_1\cr\vdots\cr
a_{n-k}\cr\end{array}\right)$ is the first column of $T$ and
$T_2=\left(\begin{array}{c}j_1\cr\vdots\cr
j_k\cr\end{array}\right)$ is the second column of $T.$ It is
enough to define the columns as sets since the entries increase
from  top  to  bottom in the columns. We denote a column by
$\langle T_i\rangle$ when we consider it as a set.

Put $\sigma_{\scriptscriptstyle T}=(i_1,j_1)\ldots (i_k,j_k)$ where
$i_1=j_1-1,$ and $i_s=\max\{d\in \langle
T_1\rangle\setminus\{i_1,\ldots,i_{s-1}\}\ |\ d<j_s\}$ for any
$s>1.$ For example, take
$$T=\begin{array}{ll}
1&4\\
2&5\\
3&7\\
6&8\\
\end{array}$$
Then $\sigma_{\scriptscriptstyle T}=(3,4)(2,5)(6,7)(1,8).$

Put $\B_T:=\B_{\sigma_T}.$ As it was shown in \cite[4.13]{Msmith}
\begin{prop}
For $T\in \bT_n^2$ one has $\overline\V_T=\overline\B_T.$
\end{prop}

 By
\cite[Remark 5.5]{M-PI} for $\sigma\in \bS_n^2(k)$ one has $\dim
\B_\sigma\leq k(n-k)$ and the equality is satisfied iff
$\sigma=\sigma_{\scriptscriptstyle T}$ for some
$T\in\bT_{(n-k,k)^*}$. That is the number of
$\B_\sigma\in\Bb_{(n-k,k)}$ of dimension $k(n-k)$ is equal to the
number of Young tableaux in $\bT_{(n-k,k)^*}$ which is ${n!\over
k!(n-k)!}{n-2k+1\over n-k+1}.$

 Let $\dim
\B_\sigma= k(n-k)$ where $\sigma=(i_1,j_1)\ldots(i_k,j_k)$ then
$\sigma=\sigma_T$ where $\langle T_2\rangle=\{j_1,\ldots, j_k\}$ and
respectively $\langle T_1\rangle=\{i\}_{i=1}^n\setminus\langle
T_2\rangle.$
\subsection{}\label{2.3}
In \cite{Mo2} the combinatorial description of $\overline\B_\sigma$
(with respect to Zariski topology) for $\sigma\in\bS_n^2$ is
provided. Let us formulate this result.

For $1\leq i<j\leq n$ consider the canonical projections
$\pi_{i,j}:\nil_n\rightarrow \nil_{j-i+1}$ acting on a matrix by
deleting the first $i-1$ columns and rows and the last $n-j$ columns
and rows. Define the rank matrix $R_u$ of $u\in\nil$ to be

$$
(R_u)_{i,j}=\left\{ \begin{array}{ll}
{\rm Rank}\,(\pi_{i,j}(u))&{\rm if}\ i< j;\\
                      0 &{\rm otherwise}.\\
                       \end{array}\right.
$$

\noindent Obviously for any $y\in\B_u$ one has $R_y=R_u$ so that we
can define $R_{\B_u}:=R_u.$

 Put $R_\sigma:=R_{\B_\sigma}=R_{N_\sigma}.$
Note that ${\rm Rank}\, (\pi_{i,j}(N_\sigma))$ is simply the number
of ones in matrix $\pi_{i,j}(N_\sigma)$ since all ones belong to
different rows and columns.

Let $\Ze^+$ be the set of non-negative integers. Put
$\bR^2_n:=\{R_\sigma\ |\ \sigma\in \bS_n^2\}.$ By \cite[3.1,
3.3]{Mo2} one has
\begin{prop} $R\in M_{n\times n}(\Ze^+)$  belongs to
$\bR^2_n$ if and only if it satisfies
\begin{itemize}
\item[(i)] $R_{i,j}=0$ if\ \ $i\geq j;$
\item[(ii)] For $i<j$ one has $R_{i+1,j}\leq R_{i,j}\leq R_{i+1,j}+1$
            and $R_{i,j-1}\leq R_{i,j}\leq R_{i,j-1}+1;$
\item[(iii)] If $R_{i,j}=R_{i+1,j}+1=R_{i,j-1}+1=R_{i+1,j-1}+1$ then
\begin{itemize}
\item[(a)] $R_{i,k}=R_{i+1,k}$ for any $k<j$ and
$R_{i,k}=R_{i+1,k}+1$ for any $k\geq j;$
\item[(b)] $R_{k,j}=R_{k,j-1}$ for any $k>i$ and
$R_{k,j}=R_{k,j-1}+1$ for any $k\leq i;$
\item[(c)] $R_{j,k}=R_{j+1,k}$ and $R_{k,i}=R_{k,i-1}$
for any $k\ :\ 1\leq k\leq n.$
\end{itemize}
\end{itemize}
\end{prop}
Note that we respectively can define $\pi_{i,j}:\bT_n\rightarrow
\bT_{j-i+1}$ by so called "jeu de taquin" (cf. \cite{Mo2}, for
example). Here we need only $\pi_{1,k}:\bT_n\rightarrow\bT_k$ where
$\pi_{1,k}(T)$ is obtained by erasing boxes with $n,n-1,\ldots,k+1.$
As it is shown in  \cite{Mo2} and one can easily see
$\pi_{i,j}\B_T=\B_{\pi_{i,j}(T)}.$

\subsection{}\label{2.4}
Recall  a partial order $\preceq$ on $\bM_n(\Real)$ defined in
\ref{1.5}. Its restriction to $\bR_n^2$ defines a partial order
$\preceq$ on $\bS_n^2$ by putting $\sigma'\preceq \sigma$ if
$R_{\sigma'}\preceq R_\sigma.$ As it is shown in \cite[3.5]{Mo2}

\begin{thm} For any $\sigma\in \bS_n^2$ one has
$$\overline\B_\sigma=\coprod\limits_{\sigma\pr\preceq\sigma}\B_{\sigma\pr}.$$
\end{thm}

\subsection{}\label{2.4a} Consider $\pi_{i,j}(\sigma)$ as an element
of $\bS_n^2$. Note that
$(R_{\pi_{i,j}(\sigma)})_{s,t}=(R_\sigma)_{\max(i,s),\min(j,t)}$ which gives us
a very easy lemma that we need in what follows.

\begin{lem} For any $\sigma,\sigma'\in\bS_n^2$ such that
$\sigma\preceq\sigma'$ and for any $i,j\ :\ 1\leq i<j\leq n$ one has
\begin{itemize}
 \item[(i)] $\pi_{i,j}(\sigma)\preceq\pi_{i,j}(\sigma')$
 \item[(ii)] $\pi_{i,j}(\sigma)\preceq\sigma$
 \end{itemize}
 \end{lem}

\subsection{}\label{2.5} As it is shown in \cite[5.6]{M-PI} there is
a unique minimal element $\sigma_o(k)$ in $\bS_n^2(k)$ (that is
$\sigma\succeq\sigma_o(k)$ for any $\sigma\in\bS_n^2(k)$). In
particular for any  $S'\subseteq\bS_n^2(k)$ one has
$$\bigcap\limits_{\sigma\in S'}\ov\B_\sigma
\bigcap\Or_{(n-k,k)^*}\ne\emptyset.\eqno{(*)}$$

To find  $\ov\B_\sigma\cap\ov\B_{\sigma'}$ for
$\sigma,\sigma'\in\bS_n^2$ we use an intersection matrix
$R_{\sigma,\sigma\pr}$ for $\sigma,\sigma\pr\in \bS_n^2$ defined in \ref{1.8}.
As it is shown in \cite[5.7]{M-PI},
\begin{thm}
For any $\sigma,\sigma\pr\in \bS_n^2$ one has
$$\overline
\B_\sigma\cap\overline\B_{\sigma\pr}=\coprod\limits_{R_{\sigma^{\prime\prime}}\preceq
R_{\sigma,\sigma'}}\B_{\sigma^{\prime\prime}}$$
%\bigcup
%\limits_{{\sigma^{\prime\prime}\in\bS_n^2(k)}\atop{R_{\sigma^{\prime\prime}}\preceq
%R_{\sigma,\sigma'}}}\ov\B_{\sigma^{\prime\prime}}.
In particular, $\overline\B_\sigma\cap\overline\B_{\sigma\pr}$ is
 irreducible if and only if $R_{\sigma,\sigma\pr}\in \bR_n^2.$ In
 that case
 $\overline\B_\sigma\cap\overline\B_{\sigma'}=\overline\B_\tau$
 where $R_\tau=R_{\sigma,\sigma\pr}.$
\end{thm}

\section{$\B_\sigma$ in terms of link patterns}
\subsection{}\label{3.1}
In this section we translate all the results of the previous
section into the language of link patterns. Recall
$\ell(P_\sigma),\ c(P_\sigma), f(P_\sigma)$ from \ref{1.4} and
$[i,j]$ from \ref{1.5}. We put $[i,j]:=\emptyset$ if $i>j.$
 For any point $m:\ 1\leq m\leq n$ put ${\bf over}_m(P_\sigma)$ to be the set of arcs over $m$, that is ${\bf over}_m(P_\sigma):=\{(i,j)\in P_\sigma\ |\
 m\in[i+1,j-1]\}.$
 For any $(i,j)\in P_\sigma$ put ${\bf over}_{(i,j)}(P_\sigma)$ to be the set of arcs over arc $(i,j)$, that is
 ${\bf over}_{(i,j)}(P_\sigma):=\{(i',j')\in
 P_\sigma\ |\  (i,j)\in[i'+1,j'-1]\}$
\begin{thm}
For $\sigma\in\bS_n^2$ one has
$$\dim
\B_\sigma=\ell(P_\sigma)\cdot(n-\ell(P_\sigma))-c(P_\sigma)-f(P_\sigma).$$
\end{thm}
\begin{proof}
Our proof is based on the recalculating  the expression of Theorem
\ref{2.1} in the language of link patterns.

 Let
$\sigma=(i_1,j_1)\ldots(i_k,j_k).$ For a given arc $(i,j)\in
P_\sigma$ we define \begin{itemize} \item{} $l_{(i,j)}(P_\sigma)$ to
be the set arcs to the left of $(i,j)$, that is
$l_{(i,j)}(P_\sigma):=\{(i',j')\in P_\sigma\ |\ (i',j')\in
[1,i-1]\};$ \item{} $r_{(i,j)}(P_\sigma)$ to be the set of arcs to
the right of $(i,j)$, that is $r_{(i,j)}(P_\sigma):=\{(i',j')\in
P_\sigma\ |\ (i',j')\in [j+1,n]\}$. \item{} ${\bf
under}_{(i,j)}(P_\sigma)$ to be the set of arcs under $(i,j)$ that
is ${\bf under}_{(i,j)}(P_\sigma):=\{(i',j')\in P_\sigma\ |\
(i',j')\in[i+1,j-1]\}$. \item{} $c^l_{(i,j)}(P_\sigma)$ to be the
set of arcs intersecting with $(i,j)$ on the left that is
$c^l_{(i,j)}(P_\sigma):=\{(i',j')\in P_\sigma\ |\  i'\in [1,i-1]\
{\rm and}\ j'\in [i+1,j-1]\}.$ \item{} $c^r_{(i,j)}(P_\sigma)$ to be
the set of arcs intersecting with $(i,j)$ on the right that is
$c^r_{(i,j)}(P_\sigma):=\{(i',j')\in P_\sigma\ |\ i'\in[i+1,j-1]\
{\rm and}\ j'\in[j+1,n]\}.$
 \item{} recall that ${\bf over}_{(i,j)}(P_\sigma):=\{(i',j')\in P_\sigma\ |\
 (i,j)\in[i'+1,j'-1]\}.$
\end{itemize}
Put $|S|$ to be cardinality of set $S.$
 Let us provide a
few identities on the summing of the  cardinalities of the sets
defined above.
\begin{itemize}
\item[(i)]
 Note that the list of six types of arcs above includes all the possible positions of an arc $(i',j')\ne(i,j)$ of $P_\sigma$
 with respect to $(i,j)$. Thus, for any $(i,j)\in P_\sigma$ one has
$$|l_{(i,j)}(P_\sigma)|+| r_{(i,j)}(P_\sigma)|+|{\bf
under}_{(i,j)}(P_\sigma)| +|{\bf over}_{(i,j)}(P_\sigma)|+
|c^l_{(i,j)}(P_\sigma)|+|c^r_{(i,j)}(P_\sigma)|=k-1.$$
 \item[(ii)] Let $\{p_s\}_{s=1}^{n-2k}$ be the set of fixed points of $P_\sigma.$
  Put $f'_{(i,j)}(P_\sigma)$ to be the number of fixed points under the arc $(i,j)$
 then obviously $\sum_{s=1}^k f'_{(i_s,j_s)}(P_\sigma)=\sum_{s=1}^{n-2k}f_{p_s}(P_\sigma)=f(P_\sigma).$ Note also that
for $(i,j)\in\sigma$ any  $p\in [i+1,j-1]$ is either a fixed
point, or some end point of $(i',j')\in{\bf
under}_{(i,j)}(P_\sigma)$ or some end point of and arc
intersecting $(i,j).$ Thus, one has
$$j-i-1=f'_{(i,j)}(P_\sigma)+2|{\bf under}_{(i,j)}(P_\sigma)|+ |c^l_{(i,j)}(P_\sigma)|+
|c^r_{(i,j)}(P_\sigma)|.$$
 \item[(iii)] Note that  $(i,j)\in c^l_{(i',j')}(P_\sigma)$ iff $(i',j')\in c^r_{(i,j)}(P_\sigma)$
  so that $$\sum\limits_{s=1}^k|c^l_{(i_s,j_s)}(P_\sigma)|=
\sum\limits_{s=1}^k|c^r_{(i_s,j_s)}(P_\sigma)|=c(P_\sigma).$$
\item[(iv)] Note also that  arc $(i',j')\in l_{(i,j)}(P_\sigma)$ iff
$(i,j)\in r_{(i',j')}(P_\sigma)$ so that
$$\sum\limits_{s=1}^k|l_{(i_s,j_s)}(P_\sigma)|=\sum\limits_{s=1}^k|r_{(i_s,j_s)}(P_\sigma)|.$$
\item[(v)] Exactly in the same way $(i',j')\in {\bf
under}_{(i,j)}(P_\sigma)$ iff $(i,j)\in{\bf
over}_{(i',j')}(P_\sigma)$ so that $$\sum\limits_{s=1}^k|{\bf
under}_{(i_s,j_s)}(P_\sigma)|=\sum\limits_{s=1}^k|{\bf
over}_{(i_s,j_s)}(P_\sigma)|.$$

\item[(vi)] Finally note that $\#\{j_p\ |\ j_p<i\}$ (from $(*)$ of
\ref{2.1}) is exactly
 $|l_{(i,j)}(P_\sigma)|$ and $\#\{i_p<i\ |\ j_p<j\}$ (from $(*)$ of \ref{2.1})
 is $|l_{(i,j)}(P_\sigma)|+|c^l_{(i,j)}(P_\sigma)|$ so that
 $$q_{(i,j)}(\sigma)=2|l_{(i,j)}(P_\sigma)|+|c^l_{(i,j)}(P_\sigma)|.$$
 \end{itemize}

In the following computations we omit the notation $(P_\sigma)$ in
all the sets defined above since we consider only $P_\sigma$ here.
Starting from Theorem \ref{2.1} we get
$$\begin{array}{l}
\dim \B_\sigma=kn-\sum\limits_{s=1}^k(j_s-i_s)-\sum\limits_{s=1}^k
q_{(i_s,j_s)}(\sigma)\quad{{\rm by\ (ii),(vi)}\atop=}\\
kn-\sum\limits_{s=1}^k(1+f'_{(i_s,j_s)}(\sigma)+2|{\bf
under}_{(i_s,j_s)}|+ |c^l_{(i_s,j_s)}|+|c^r_{(i_s,j_s)}|)\\
\quad-\sum\limits_{s=1}^k(2|l_{(i_s,j_s)}|+|c^l_{(i_s,j_s)}|)\quad {{{\rm by\ (iii)\ and\ def.\ }f(\sigma)  }\atop =}\\
 nk-k-f(\sigma)-3c(P_\sigma)-\sum\limits_{s=1}^k(2|{\bf
under}_{(i_s,j_s)}|+2|l_{(i_s,j_s)}|)\quad{{\rm by\ (iv),(v)}\atop =}\\
nk-k-f(\sigma)-3c(P_\sigma)-\sum\limits_{s=1}^k(|{\bf
under}_{(i_s,j_s)}|+|{\bf
over}_{(i_s,j_s)}|+|l_{(i_s,j_s)}|+|r_{(i_s,j_s)}|)\quad
{{\rm by\ (i)}\atop = }\\
nk-k-f(\sigma)-3c(P_\sigma)-\sum\limits_{s=1}^k(k-1-(|c^l_{(i_s,j_s)}|+
|c^r_{(i_s,j_s)}|))\quad{{\rm by\ (iii)}\atop
=}\\
nk-k-f(\sigma)-3c(P_\sigma)-k(k-1)+2c(P_\sigma)=\\
k(n-k)-f(P_\sigma)-c(P_\sigma).\\
 \end{array}
$$
\end{proof}

\subsection{}\label{3.2}
Let us translate the definitions of \ref{2.3} into the language of
link patterns.

 For $i,j\ :\  1\leq i<j\leq n$ we define
$\pi_{i,j}(P_\sigma)$ to be a link pattern from $i$ to $j$. It
consists only of arcs $(i',j')\in P_\sigma$ such that
$(i',j')\in[i,j].$ On the other hand, there is 1 in place $(i',j')$
of matrix $N_\sigma$ where $i'\geq i$ and $j'\leq j$ iff there
exists $(i',j')\in\sigma$. Hence $\ell(\pi_{i,j}(P_\sigma))$ is
equal to number of ones in $\pi_{i,j}(N_\sigma)$ which in turn is
${\rm Rank}(\pi_{i,j}(N_\sigma))$ as it was noted in \ref{2.3}. Thus,
$$
 (R_\sigma)_{i,j}=\left\{
\begin{array}{ll}
\ell(\pi_{i,j}(P_\sigma))& {\rm if}\ i<j;\\
0&{\rm otherwise;}\\
\end{array}\right.$$
exactly as we have defined in \ref{1.5}. Together with \ref{2.4}
this provides the result announced in \ref{1.5}.

\subsection{}\label{3.3} Now we can describe the cover of $\sigma$
for our partial order in terms of link patterns. We start with the
description of $\N(\sigma)$ defined in \ref{1.6}.

Recall from \ref{1.6} that  $(i,j)\in P_\sigma$ is called external
if ${\bf over}_{(i,j)}(P_\sigma)=\emptyset.$ Let $E(\sigma)$ be the
set of external arcs in $P_\sigma$.

For example, let $\sigma=(1,3)(2,7)(4,5)$ then $P_\sigma$ is
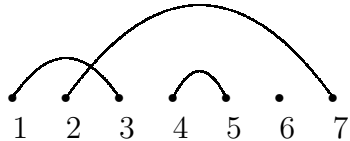
\begin{center}
\begin{picture}(100,80)
\multiput(10,40)(20,0){7}%
{\circle*{3}}
 \put(10,25){1}
 \put(30,25){2}
 \put(50,25){3}
 \put(70,25){4}
 \put(90,25){5}
 \put(110,25){6}
 \put(130,25){7}
 \qbezier(10,40)(30,70)(50,40)
 \qbezier(30,40)(80,110)(130,40)
 \qbezier(70,40)(80,60)(90,40)
\end{picture}
\end{center}
and $E(\sigma)=\{(1,3),(2,7)\}.$

For $(i,j)\in\sigma$ put $\sigma_{(i,j)}^-$ to be an involution
obtained from $\sigma$ by omitting $(i,j)$. Respectively
$P_{\sigma_{(i,j)}^-}$ is obtained from $P_\sigma$ by erasing arc
$(i,j).$

 By \cite[3.7, 3.9]{Mo2} one has
\begin{prop}
 $\N(\sigma)=\{\sigma_{(i,j)}^-\ |\ (i,j)\in E(\sigma)\}$
\end{prop}

\subsection{}\label{3.4}
$\D(\sigma)$ defined in \ref{1.6} is more involved. Its description
given in \cite[3.11-3.14]{Mo2} is rather complicated. However it is
much more short and elegant in the language of link patterns.
$\D(\sigma)$ contains of three types of elements which we describe
in three following subsections. The main idea which is very clear on
the level of link patterns is that the elements of $\D(\sigma)$ are
obtained either by moving a left (resp. right) end point of an arc
to the nearest fixed point to the left (resp. to the right) of it,
so that we get a new fixed point under the arc, or by crossing two
arcs. If $P_{\sigma'}$ is obtained from $P_\sigma$ in such way, then
$\B_{\sigma'}$ must be of codimension 1 in $\ov\B_\sigma.$ A more
delicate point is that all the elements of $\D(\sigma)$ are obtained
in such way. This is exactly the content of Proposition 3.15 of
\cite{Mo2}.

We begin with constructing of a link pattern obtained from a given
one by moving one end point of an arc to the nearest fixed point.

 Given
$\sigma=(i_1,j_1)\ldots(i_k,j_k)$ we put $\langle
\sigma\rangle:=\langle P_\sigma\rangle:=\{i_1,j_1,\ldots i_k,j_k\}$
to be the set of end points of $P_\sigma$.

For $i\in \langle \sigma\rangle$ and $f\not\in \langle
\sigma\rangle$ put $\sigma_{i\rightarrow f}$ to be an involution
obtained from $\sigma$ by changing an arc $((i,j))$ to $((f,j))$.
We write the arcs in double round brackets since the definition is not connected to the ordering of the ends of an arc.

For $(i,j)\in \sigma$ if there exist fixed points to the left of arc
$(i,j)$ let $m$ be the largest among them, that is $m=\max\{p<i\ |\
p\not\in \langle \sigma\rangle\}.$ If in addition ${\bf
over}_{(i,j)}(P_\sigma)\cap [m+1,n]=\emptyset$ (or in other words
${\bf over}_{(i,j)}(P_\sigma)\subseteq {\bf over}_m(P_\sigma)$) we
put $\sigma_{\curvearrowleft (i,j)}:=\sigma_{i\rightarrow m}.$
Otherwise put $\sigma_{\curvearrowleft (i,j)}=\emptyset.$

Let $\sigma_{\curvearrowleft(i,j)}=\sigma_{i\rightarrow m}$. Note
that by \ref{3.2} we get
$$ (R_{\sigma_{\curvearrowleft(i,j)}})_{s,t}=\left\{\begin{array}{ll}
(R_\sigma)_{s,t}-1 &{\rm if}\ m< s\leq i\ {\rm and}\ t\geq j;\\
(R_\sigma)_{s,t}&{\rm otherwise.}\\
\end{array}\right.$$ Thus, if
$\sigma_{\curvearrowleft(i,j)}\ne\emptyset$ then
$\sigma_{\curvearrowleft(i,j)}\prec\sigma.$

 Exactly in the same way, for $(i,j)\in \sigma$ if there exist fixed points to the right of $(i,j)$
 let $m$ be the minimal among them, that is
$m=\min\{p>j\ |\ p\not\in \langle \sigma\rangle\}.$ If in addition
${\bf over}_{(i,j)}(P_\sigma)\cap [1,m-1]=\emptyset$ (that is again
${\bf over}_{(i,j)}(P_\sigma)\subseteq {\bf over}_m(P_\sigma)$)
 we put
$\sigma_{(i,j)\curvearrowright}:=\sigma_{j\rightarrow m}$. Otherwise
put $\sigma_{ (i,j)\curvearrowright}=\emptyset.$

And again, exactly as in case $\sigma_{\curvearrowleft (i,j)}$ if
$\sigma_{(i,j)\curvearrowright}\ne\emptyset$ then
$\sigma_{(i,j)\curvearrowright}\prec\sigma.$

 For example, let $\sigma=(1,6)(3,5)(4,7)\in\bS_8^2.$
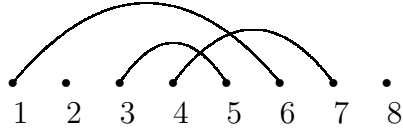
\begin{center}
\begin{picture}(100,80)
\multiput(-40,40)(20,0){8}%
{\circle*{3}}
 \put(-40,25){1}
 \put(-20,25){2}
 \put(0,25){3}
 \put(20,25){4}
 \put(40,25){5}
 \put(60,25){6}
 \put(80,25){7}
 \put(100,25){8}

 \qbezier(-40,40)(10,100)(60,40)
 \qbezier(0,40)(20,70)(40,40)
 \qbezier(20,40)(50,80)(80,40)
\end{picture}
\end{center}
Then $\sigma_{\curvearrowleft (1,6)}=\emptyset$ and
$\sigma_{(1,6)\curvearrowright}=(1,8)(3,5)(4,7)$ so that
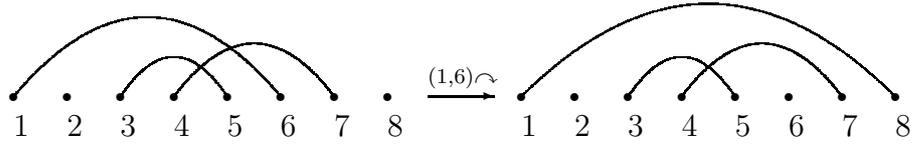
\begin{center}
\begin{picture}(100,80)
\multiput(-100,40)(20,0){8}%
{\circle*{3}}
 \put(-100,25){1}
 \put(-80,25){2}
 \put(-60,25){3}
 \put(-40,25){4}
 \put(-20,25){5}
 \put(0,25){6}
 \put(20,25){7}
 \put(40,25){8}

 \qbezier(-100,40)(-50,100)(0,40)
 \qbezier(-60,40)(-40,70)(-20,40)
 \qbezier(-40,40)(-10,80)(20,40)
 \put(55,45){$\scriptstyle{(1,6)\curvearrowright}$}
 \put(55,40){\vector(1,0){25}}
 \multiput(90,40)(20,0){8}%
 {\circle*{3}}
\put(90,25){1}
 \put(110,25){2}
 \put(130,25){3}
 \put(150,25){4}
 \put(170,25){5}
 \put(190,25){6}
 \put(210,25){7}
 \put(230,25){8}

 \qbezier(90,40)(160,110)(230,40)
 \qbezier(130,40)(150,70)(170,40)
 \qbezier(150,40)(180,80)(210,40)

\end{picture}
\end{center}
$\sigma_{\curvearrowleft (3,5)}=(1,6)(2,5)(4,7)$ and
$\sigma_{(3,5)\curvearrowright}=\emptyset$ so that
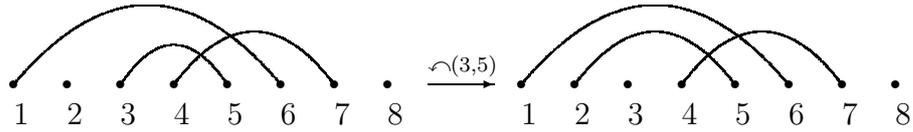
\begin{center}
\begin{picture}(100,80)
\multiput(-100,40)(20,0){8}%
{\circle*{3}}
 \put(-100,25){1}
 \put(-80,25){2}
 \put(-60,25){3}
 \put(-40,25){4}
 \put(-20,25){5}
 \put(0,25){6}
 \put(20,25){7}
 \put(40,25){8}

 \qbezier(-100,40)(-50,100)(0,40)
 \qbezier(-60,40)(-40,70)(-20,40)
 \qbezier(-40,40)(-10,80)(20,40)
 \put(55,45){$\scriptstyle{\curvearrowleft (3,5)}$}
 \put(55,40){\vector(1,0){25}}
 \multiput(90,40)(20,0){8}%
 {\circle*{3}}
\put(90,25){1}
 \put(110,25){2}
 \put(130,25){3}
 \put(150,25){4}
 \put(170,25){5}
 \put(190,25){6}
 \put(210,25){7}
 \put(230,25){8}

 \qbezier(90,40)(140,100)(190,40)
 \qbezier(110,40)(140,80)(170,40)
 \qbezier(150,40)(180,80)(210,40)

\end{picture}
\end{center}

Finally, $\sigma_{\curvearrowleft (4,7)}=(1,6)(3,5)(2,7)$ and
$\sigma_{(4,7)\curvearrowright}=(1,6)(3,5)(4,8)$ so that
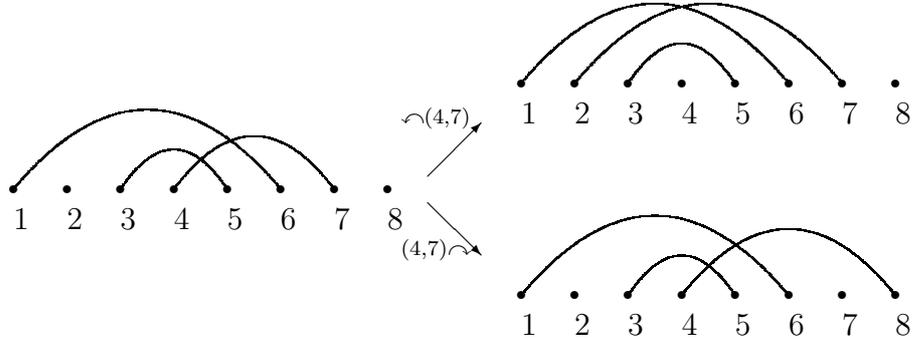
\begin{center}
\begin{picture}(100,150)
\multiput(-100,60)(20,0){8}%
{\circle*{3}}
 \put(-100,45){1}
 \put(-80,45){2}
 \put(-60,45){3}
 \put(-40,45){4}
 \put(-20,45){5}
 \put(0,45){6}
 \put(20,45){7}
 \put(40,45){8}

 \qbezier(-100,60)(-50,120)(0,60)
 \qbezier(-60,60)(-40,90)(-20,60)
 \qbezier(-40,60)(-10,100)(20,60)
 %\put(60,40){${\curvearrowleft (3,5)}\atop{\ \longrightarrow\ }$}
\put(55,65){\vector(1,1){20}}
 \put (55,55){\vector(1,-1){20}}
\put(45,35){$\scriptstyle{ (4,7)\curvearrowright}$}
\put(45,85){$\scriptstyle{\curvearrowleft (4,7)}$}

 \multiput(90,100)(20,0){8}%
 {\circle*{3}}
 \put(90,85){1}
 \put(110,85){2}
 \put(130,85){3}
 \put(150,85){4}
 \put(170,85){5}
 \put(190,85){6}
 \put(210,85){7}
 \put(230,85){8}

 \qbezier(90,100)(140,160)(190,100)
 \qbezier(130,100)(150,130)(170,100)
 \qbezier(110,100)(160,160)(210,100)
\multiput(90,20)(20,0){8}%
 {\circle*{3}}
 \put(90,5){1}
 \put(110,5){2}
 \put(130,5){3}
 \put(150,5){4}
 \put(170,5){5}
 \put(190,5){6}
 \put(210,5){7}
 \put(230,5){8}

 \qbezier(90,20)(140,80)(190,20)
 \qbezier(130,20)(150,50)(170,20)
 \qbezier(150,20)(190,70)(230,20)

\end{picture}
\end{center}

As an easy corollary of Theorem \ref{3.1} we get
\begin{prop} Let $\sigma'$ be either $\sigma_{\curvearrowleft (i,j)}$
or $\sigma_{ (i,j)\curvearrowright}.$ Then
$\codim_{\ov\B_\sigma}\B_{\sigma'}=1.$
\end{prop}
\begin{proof}
By our note above $\B_{\sigma'}\subset\ov\B_\sigma.$

Let us compute $\dim \B_{\sigma'}.$ It is enough to compute it  in
one of the cases since one case can be obtained from the other by
mirroring $P_\sigma$ around the point ${1+n}\over 2$.

 For example, let us
consider case $\sigma'=\sigma_{\curvearrowleft (i,j)}.$ Let
 $\sigma'=\sigma_{i\rightarrow m}.$ Recall notation from the
 proof of Theorem \ref{3.1}. Note that
\begin{itemize}
\item{} $f'_{(m,j)}(P_{\sigma'})=f'_{(i,j)}(P_\sigma)+1$ and since
${\bf over}_{(i,j)}(P_\sigma)\subset {\bf over}_m(P_\sigma)$ one has
$c^r_{(m,j)}(P_{\sigma'})=c^r_{(i,j)}(P_\sigma).$

\item{} Again, since ${\bf over}_{(i,j)}(P_\sigma)\subset {\bf over}_m(P_\sigma)$
one has  ${\bf over}_{(m,j)}(P_{\sigma'})={\bf
over}_{(i,j)}(P_\sigma)$ so that for any $(i',j')\in {\bf
over}_{(i,j)}(P_\sigma)$ one has
$f'_{(i',j')}(P_{\sigma'})=f'_{(i',j')}(P_\sigma)$ and
$c^r_{(i',j')}(P_{\sigma'})=c^r_{(i',j')}(P_\sigma).$
\item{} Exactly in the same way for $(i',j')\in P_\sigma$ such that
 $(i',j')\not\in{\bf over}_m(P_\sigma)$ one has
 $f'_{(i',j')}(P_{\sigma'})=f'_{(i',j')}(P_\sigma)$ and
$c^r_{(i',j')}(P_{\sigma'}=c^r_{(i',j')}(P_\sigma).$
\item{} If $(i',j')\in\left({\bf over}_m(P_\sigma)\setminus {\bf over}_{(i,j)}(P_\sigma)\right)\cap
c^l_{(i,j)}(P_\sigma)$ (that is  $i'<m$ and $i<j'<j$) then again
$f'_{(i',j')}(P_{\sigma'})=f'_{(i',j')}(P_\sigma)$ and
$c^r_{(i',j')}(P_{\sigma'}=c^r_{(i',j')}(P_\sigma).$
\item{} Finally, if $(i',j')\in\left({\bf over}_m(P_\sigma)\setminus {\bf over}_{(i,j)}(P_\sigma)
\right)\cap l_{(i,j)}(P_\sigma)$ (that is $i'<m$ and $m<j'<i$) then
$f'_{(i',j')}(P_{\sigma'})=f'_{(i',j')}(P_\sigma)-1$ and
$c^r_{(i',j')}(P_{\sigma'})=c^r_{(i',j')}(P_\sigma)+1.$
\end{itemize}
Set $t:=|({\bf over}_m(P_\sigma)\setminus {\bf
over}_{(i,j)}(P_\sigma))\cap l_{(i,j)}(P_\sigma)|.$ Summarizing the
calculations above we get
$$\begin{array}{ll}
f(P_{\sigma'})&=\sum_{s=1}^k
f'_{(i_s,j_s)}(P_{\sigma'})=f(P_\sigma)+1-t;\\
c(P_{\sigma'})&=\sum_{s=1}^k
c^r_{(i_s,j_s)}(P_{\sigma'})=c(P_\sigma)+t;\\
\end{array}$$
so that $\dim \B_{\sigma'}=k(n-k)-f(P_{\sigma'})-c(P_{\sigma'})=\dim
\B_{\sigma}-1.$
\end{proof}

\subsection{}\label{3.5} Let us again consider
$\sigma=(i_1,j_1)\ldots(i_k,j_k).$ Take  $i,j\in\langle
\sigma\rangle$ such that $i<j.$  Assume that they belong to
different pairs, that is $i\in (i_s,j_s)$ and $j\in(i_t,j_t)$ where
$s\ne t.$ Put $\sigma_{i\leftrightarrows j}$ to be an involution
obtained by interchanging places of $i$ and $j$ in the pairs. Assume
that $(i_s,j_s)=((i,p))$ and $(i_t,j_t)=((j,q))$ then
$\sigma_{i\leftrightarrows j}=((i,q))((j,p))\ldots$ where by dots we
denote all the pairs of $\sigma$ but $(i_s,j_s),(i_t,j_t).$ Note
that we cannot say anything about ordering $((i,q))$ and $((j,p))$
so we write them in double round brackets.

The link pattern $P_{\sigma'}$ described below is obtained from
$P_\sigma$ by crossing some arc with an arc to the left of it.

Let $(i,j)\in\sigma.$ Translating
 \cite[3.13]{Mo2} we put
$$L_{(i,j)}(\sigma):=\{(i_s,j_s)\in\sigma\ |\ j_s<i\ {\rm and}\
\{p\}_{p=j_s+1}^{i-1}\subset\langle \pi_{i_s,j}(P_\sigma)\rangle\}$$
In other words $(i_s,j_s)\in L_{(i,j)}(\sigma)$ if $(i_s,j_s)\in
l_{(i,j)}(P_\sigma)$ and there are no fixed points between $j_s$ and
$i$ in $\pi_{i_s,j}(P_\sigma).$ Note that $L_{(i,j)}(\sigma)$ may
contain a few elements. Put
$$S_{(i,j)\looparrowright}(\sigma):=\left\{\begin{array}{ll}
\{\sigma_{j_s\leftrightarrows i}\ |\ (i_s,j_s)\in
L_{(i,j)}(\sigma)\}& {\rm
if}\ L_{(i,j)}(\sigma)\ne \emptyset\\
\emptyset& {\rm otherwise}\\
\end{array}\right.$$

Note that for $(i_s,j_s)\in L_{(i,j)}(\sigma)$ one has
$$(R_{\sigma_{j_s\leftrightarrows i}})_{s,t}=\left\{\begin{array}{ll}
(R_\sigma)_{s,t}-1 &{\rm if}\ i\leq s<j_s\ {\rm and}\ t\geq j\
{\rm or}\ j_s\leq t<i\ {\rm and}\ s\leq i_s ;\\
(R_\sigma)_{s,t}&{\rm otherwise.}\\
\end{array}\right.$$ Thus, for any $(i_s,j_s)\in L_{(i,j)}(\sigma)$
one has $\sigma_{j_s\leftrightarrows i}\prec\sigma.$

For example, let $\sigma=(1,5)(2,4)(3,6)(7,9)(10,11)\in\bS_{11}^2.$
Its link pattern is
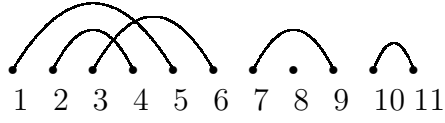
\begin{center}
\begin{picture}(100,80)
\multiput(-40,40)(15,0){11}%
{\circle*{3}}
 \put(-40,25){1}
 \put(-25,25){2}
 \put(-10,25){3}
 \put(5,25){4}
 \put(20,25){5}
 \put(35,25){6}
 \put(50,25){7}
 \put(65,25){8}
\put(80,25){9}
 \put(95,25){10}
 \put(110,25){11}

 \qbezier(-40,40)(-10,90)(20,40)
 \qbezier(-25,40)(-10,70)(5,40)
\qbezier(-10,40)(12.5,80)(35,40)
 \qbezier(50,40)(65,70)(80,40)
\qbezier(95,40)(102.5,60)(110,40)
\end{picture}
\end{center}
One has
$L_{(1,5)}(\sigma)=L_{(2,4)}(\sigma)=L_{(3,6)}(\sigma)=\emptyset,$
$L_{(7,9)}(\sigma)=\{(1,5),\ (3,6)\}$,  so that

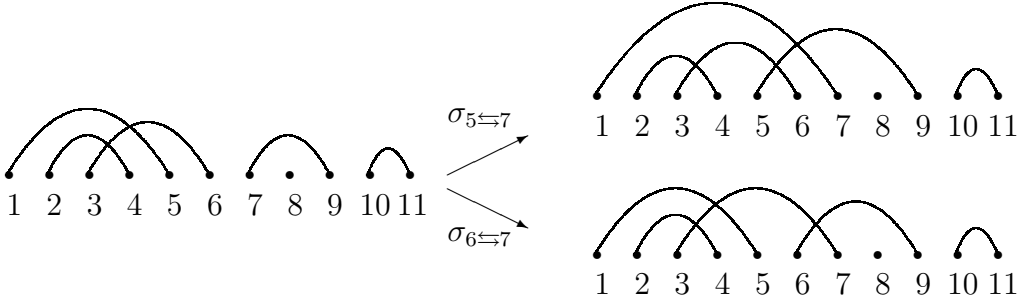
\begin{center}
\begin{picture}(100,140)
\multiput(-139,60)(15,0){11}%
{\circle*{3}}
 \put(-140,45){1}
 \put(-125,45){2}
 \put(-110,45){3}
 \put(-95,45){4}
 \put(-80,45){5}
 \put(-65,45){6}
 \put(-50,45){7}
 \put(-35,45){8}
\put(-20,45){9}
 \put(-8,45){10}
 \put(6,45){11}

 \qbezier(-140,60)(-110,110)(-79,60)
 \qbezier(-125,60)(-110,90)(-94,60)
\qbezier(-110,60)(-87.5,100)(-64,60)
 \qbezier(-50,60)(-35,90)(-19,60)
\qbezier(-5,60)(2.5,80)(11,60)

\put(25,80){$\sigma_{5\leftrightarrows 7}$}
 \put(25,60){\vector(2,1){30}}

\multiput(81,90)(15,0){11}%
{\circle*{3}}
 \put(80,75){1}
 \put(95,75){2}
 \put(110,75){3}
 \put(125,75){4}
 \put(140,75){5}
 \put(155,75){6}
 \put(170,75){7}
 \put(185,75){8}
\put(200,75){9}
 \put(212,75){10}
 \put(227,75){11}

 \qbezier(80,90)(125,160)(171,90)
 \qbezier(95,90)(110,120)(126,90)
\qbezier(140,90)(170,140)(201,90)
 \qbezier(110,90)(132.5,130)(156,90)
\qbezier(215,90)(222.5,110)(231,90)

\put(25,35){$\sigma_{6\leftrightarrows 7}$}
 \put(25,55){\vector(2,-1){30}}

 \multiput(81,30)(15,0){11}%
{\circle*{3}}
 \put(80,15){1}
 \put(95,15){2}
 \put(110,15){3}
 \put(125,15){4}
 \put(140,15){5}
 \put(155,15){6}
 \put(170,15){7}
 \put(185,15){8}
\put(200,15){9}
 \put(212,15){10}
 \put(227,15){11}

 \qbezier(80,30)(110,80)(141,30)
 \qbezier(95,30)(110,60)(126,30)
\qbezier(110,30)(140,80)(171,30)
 \qbezier(155,30)(177.5,70)(201,30)
\qbezier(215,30)(222.5,50)(231,30)

\end{picture}
\end{center}
and $L_{(10,11)}(\sigma)=\{(7,9)\}$ so that
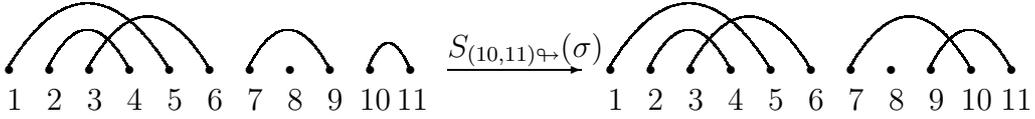
\begin{center}
\begin{picture}(100,80)
\multiput(-139,40)(15,0){11}%
{\circle*{3}}
 \put(-140,25){1}
 \put(-125,25){2}
 \put(-110,25){3}
 \put(-95,25){4}
 \put(-80,25){5}
 \put(-65,25){6}
 \put(-50,25){7}
 \put(-35,25){8}
\put(-20,25){9}
 \put(-8,25){10}
 \put(6,25){11}

 \qbezier(-140,40)(-110,90)(-79,40)
 \qbezier(-125,40)(-110,70)(-94,40)
\qbezier(-110,40)(-87.5,80)(-64,40)
 \qbezier(-50,40)(-35,70)(-19,40)
\qbezier(-5,40)(2.5,60)(11,40)

\put(25,45){$S_{(10,11)\looparrowright}(\sigma)$}
 \put(25,40){\vector(1,0){50}}

\multiput(86,40)(15,0){11}%
{\circle*{3}}
 \put(85,25){1}
 \put(100,25){2}
 \put(115,25){3}
 \put(130,25){4}
 \put(145,25){5}
 \put(160,25){6}
 \put(175,25){7}
 \put(190,25){8}
\put(205,25){9}
 \put(217,25){10}
 \put(232,25){11}

 \qbezier(85,40)(115,90)(146,40)
 \qbezier(100,40)(115,70)(131,40)
\qbezier(115,40)(137.5,80)(161,40)
 \qbezier(175,40)(197.5,80)(221,40)
\qbezier(205,40)(220,70)(236,40)
\end{picture}
\end{center}
Again, as an easy corollary of Theorem \ref{3.1} we get
\begin{prop} If $L_{(i,j)}(\sigma)\ne\emptyset$ then
$\codim_{\ov\B_\sigma}\B_{\sigma'}=1$ for any $\sigma'\in
S_{(i,j)\looparrowright}(\sigma)$
\end{prop}
\begin{proof}
By the note above $\B_{\sigma'}\st\ov\B_\sigma.$ Also, by the
definition of $L_{(i,j)}(\sigma)$ $P_{\sigma'}$ differs from
$P_\sigma$ by one additional cross obtained from changing pairs
$(i_s,j_s)(i,j)$ to $(i_s,i)(j_s,j)$ as we show below.

 Indeed, since there are no
fixed points between $j_s$ and $i$ in $\pi_{i_s,j}(P_\sigma)$ we get
that in particular there are no fixed points between $j_s$ and $i$
in $P_\sigma$ so that $f(P_\sigma)=f(P_{\sigma'}).$ Also, the
absence of fixed points between $j_s$ and $i$ in
$\pi_{i_s,j}(P_\sigma)$ means that ${\bf
over}_{(i_s,j_s)}(P_\sigma)\cap l_{(i,j)}(P_\sigma)=\emptyset$ and
${\bf over}_{(i,j)}(P_\sigma)\cap r_{(i_s,j_s)}(P_\sigma)=\emptyset$
 so that $ c(P_{\sigma'})=c(P_\sigma)+1$. Thus, by Theorem \ref{3.1}
$\dim \B_{\sigma'}=\dim\B_\sigma-1.$
\end{proof}

\subsection{}\label{3.6} The last type of elements in $\D(\sigma)$
is obtained by crossing two concentric arcs. It is defined as
follows. For $(i,j)\in\sigma$ recall ${\bf over}_{(i,j)}(P_\sigma)$
and ${\bf under}_{(i,j)}(P_\sigma)$ from \ref{3.1} and put
$$Ov_{(i,j)}(\sigma):=\{(i_s,j_s)\in {\bf over}_{(i,j)}(P_\sigma)\
|\
 {\bf under}_{(i_s,j_s)}(P_\sigma)=\{(i,j)\}\cup{\bf under}_{(i,j)}(P_\sigma) \}$$
 In other words $(i_s,j_s)\in Ov_{(i,j)}(\sigma)$ if $(i_s,j_s)\in
{\bf over}_{(i,j)}(P_\sigma)$ and there are no arcs between
$(i_s,j_s)$ and $(i,j)$ in $P_\sigma.$ Note that
$Ov_{(i,j)}(\sigma)$ may contain a few elements. Put
$$S_{(i,j)\upharpoonleft\downharpoonright}(\sigma):=\left\{\begin{array}{ll}
\{\sigma_{i_s\leftrightarrows i}\ |\ (i_s,j_s)\in
Ov_{(i,j)}(\sigma)\} & {\rm if}\ Ov_{(i,j)}(\sigma)\ne\emptyset;\\
\emptyset &{\rm otherwise};\\
\end{array}\right.
$$
Note that $\sigma_{i_s\leftrightarrows
i_t}=\sigma_{j_s\leftrightarrows j_t}$. Note also that for
$(i_p,j_p)\in Ov_{(i,j)}(\sigma)$ one has
$$(R_{\sigma_{i_p\leftrightarrows i}})_{s,t}=\left\{\begin{array}{ll}
(R_\sigma)_{s,t}-1 &{\rm if}\ i_p< s\leq i\ {\rm and}\ j\leq t<
j_s;\\
(R_\sigma)_{s,t}&{\rm otherwise.}\\
\end{array}\right.$$
 Thus, for any $(i_p,j_p)\in Ov_{(i,j)}(\sigma)$ one has
$\sigma_{i_p\leftrightarrows i}\prec\sigma.$

 For example, let
$\sigma=(1,11)(2,6)(3,9)(4,5)\in\bS_{11}^2.$ Its link pattern is
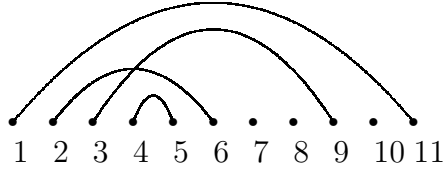
\begin{center}
\begin{picture}(100,100)
\multiput(-40,40)(15,0){11}%
{\circle*{3}}
 \put(-40,25){1}
 \put(-25,25){2}
 \put(-10,25){3}
 \put(5,25){4}
 \put(20,25){5}
 \put(35,25){6}
 \put(50,25){7}
 \put(65,25){8}
\put(80,25){9}
 \put(95,25){10}
 \put(110,25){11}

 \qbezier(-40,40)(35,130)(110,40)
 \qbezier(-25,40)(5,80)(35,40)
\qbezier(-10,40)(35,110)(80,40)
 \qbezier(5,40)(12.5,60)(20,40)
\end{picture}
\end{center}
One has $Ov_{(1,11)}(\sigma)=\emptyset,$
$Ov_{(2,6)}(\sigma)=\{(1,11)\}$ so that
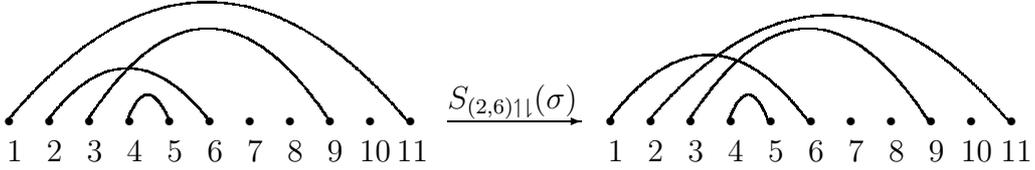
\begin{center}
\begin{picture}(100,100)
\multiput(-139,40)(15,0){11}%
{\circle*{3}}
 \put(-140,25){1}
 \put(-125,25){2}
 \put(-110,25){3}
 \put(-95,25){4}
 \put(-80,25){5}
 \put(-65,25){6}
 \put(-50,25){7}
 \put(-35,25){8}
\put(-20,25){9}
 \put(-8,25){10}
 \put(6,25){11}

 \qbezier(-140,40)(-65,130)(10,40)
 \qbezier(-125,40)(-95,80)(-64,40)
\qbezier(-110,40)(-65,110)(-19,40)
 \qbezier(-95,40)(-87.5,60)(-79,40)

\put(25,45){$S_{(2,6)\upharpoonleft\downharpoonright}(\sigma)$}
 \put(25,40){\vector(1,0){50}}

\multiput(86,40)(15,0){11}%
{\circle*{3}}
 \put(85,25){1}
 \put(100,25){2}
 \put(115,25){3}
 \put(130,25){4}
 \put(145,25){5}
 \put(160,25){6}
 \put(175,25){7}
 \put(190,25){8}
\put(205,25){9}
 \put(217,25){10}
 \put(232,25){11}

 \qbezier(85,40)(122.5,90)(160,40)
 \qbezier(100,40)(167.5,120)(235,40)
\qbezier(115,40)(160,110)(205,40)
 \qbezier(130,40)(137.5,60)(145,40)

\end{picture}
\end{center}
and $Ov_{(3,9)}(\sigma)=\{(1,11)\}$ so that
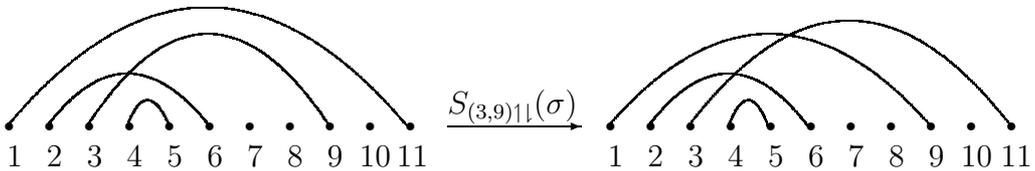
\begin{center}
\begin{picture}(100,100)
\multiput(-139,40)(15,0){11}%
{\circle*{3}}
 \put(-140,25){1}
 \put(-125,25){2}
 \put(-110,25){3}
 \put(-95,25){4}
 \put(-80,25){5}
 \put(-65,25){6}
 \put(-50,25){7}
 \put(-35,25){8}
\put(-20,25){9}
 \put(-8,25){10}
 \put(6,25){11}

 \qbezier(-140,40)(-65,130)(10,40)
 \qbezier(-125,40)(-95,80)(-64,40)
\qbezier(-110,40)(-65,110)(-19,40)
 \qbezier(-95,40)(-87.5,60)(-79,40)

\put(25,45){$S_{(3,9)\upharpoonleft\downharpoonright}(\sigma)$}
 \put(25,40){\vector(1,0){50}}

\multiput(86,40)(15,0){11}%
{\circle*{3}}
 \put(85,25){1}
 \put(100,25){2}
 \put(115,25){3}
 \put(130,25){4}
 \put(145,25){5}
 \put(160,25){6}
 \put(175,25){7}
 \put(190,25){8}
\put(205,25){9}
 \put(217,25){10}
 \put(232,25){11}

 \qbezier(85,40)(145,110)(206,40)
 \qbezier(100,40)(130,80)(160,40)
\qbezier(115,40)(175,120)(236,40)
 \qbezier(130,40)(137.5,60)(145,40)

\end{picture}
\end{center}
Finally $Ov_{(4,5)}(\sigma)=\{(2,6),\ (3,9)\}$ so that
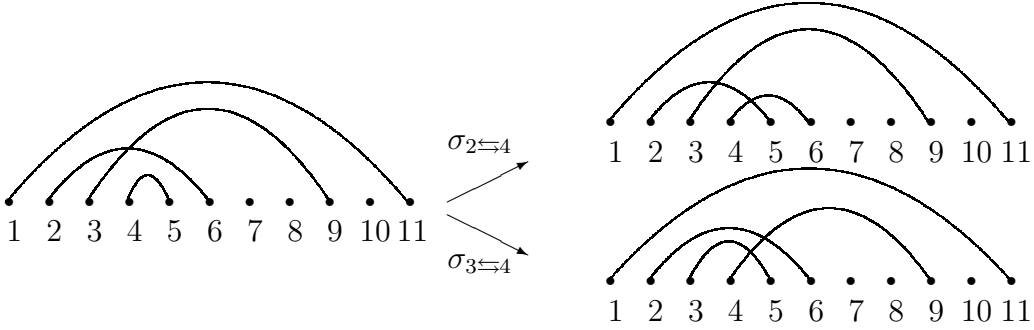
\begin{center}
\begin{picture}(100,140)
\multiput(-139,60)(15,0){11}%
{\circle*{3}}
 \put(-140,45){1}
 \put(-125,45){2}
 \put(-110,45){3}
 \put(-95,45){4}
 \put(-80,45){5}
 \put(-65,45){6}
 \put(-50,45){7}
 \put(-35,45){8}
\put(-20,45){9}
 \put(-8,45){10}
 \put(6,45){11}

 \qbezier(-140,60)(-65,150)(10,60)
 \qbezier(-125,60)(-95,100)(-64,60)
\qbezier(-110,60)(-65,130)(-19,60)
 \qbezier(-95,60)(-87.5,80)(-79,60)

\put(25,80){$\sigma_{2\leftrightarrows 4}$}
 \put(25,60){\vector(2,1){30}}

\multiput(86,90)(15,0){11}%
{\circle*{3}}
 \put(85,75){1}
 \put(100,75){2}
 \put(115,75){3}
 \put(130,75){4}
 \put(145,75){5}
 \put(160,75){6}
 \put(175,75){7}
 \put(190,75){8}
\put(205,75){9}
 \put(217,75){10}
 \put(232,75){11}

 \qbezier(85,90)(160,180)(235,90)
 \qbezier(100,90)(122.5,120)(146,90)
\qbezier(115,90)(160,160)(205,90)
 \qbezier(130,90)(145,110)(160,90)

\put(25,35){$\sigma_{3\leftrightarrows 4}$}
 \put(25,55){\vector(2,-1){30}}

 \multiput(86,30)(15,0){11}%
{\circle*{3}}
 \put(85,15){1}
 \put(100,15){2}
 \put(115,15){3}
 \put(130,15){4}
 \put(145,15){5}
 \put(160,15){6}
 \put(175,15){7}
 \put(190,15){8}
\put(205,15){9}
 \put(217,15){10}
 \put(232,15){11}

 \qbezier(85,30)(160,115)(236,30)
 \qbezier(100,30)(130,70)(161,30)
\qbezier(115,30)(130,60)(146,30)
 \qbezier(130,30)(167.5,85)(206,30)

\end{picture}
\end{center}
And again, exactly as in \ref{3.5} we get
\begin{prop} If $Ov_{(i,j)}(\sigma)\ne\emptyset$ then
$\codim_{\ov\B_\sigma}\B_{\sigma'}=1$ for any $\sigma'\in
S_{(i,j)\upharpoonleft\downharpoonright}(\sigma)$
\end{prop}
\begin{proof}
Exactly as in \ref{3.5} we get that
$\B_{\sigma'}\subset\ov\B_\sigma$ and exactly as in \ref{3.5} we get
that $f(P_{\sigma'})=f(P_\sigma)$ and $c(P_{\sigma'})=c(P_\sigma)+1$
which together provide us the result.
\end{proof}

\subsection{}\label{3.7} By Proposition 3.15 from \cite{Mo2} and
Propositions \ref{3.4}, \ref{3.5}, \ref{3.6} we get \begin{thm} Let
$\sigma=(i_1,j_1)\ldots(i_k,j_k)\in\bS_n^2(k)$ then \begin{itemize}
\item[(i)]
$\D(\sigma)=\coprod\limits_{s=1}^k\sigma_{(i_s,j_s)\curvearrowright}
\bigsqcup\coprod\limits_{s=1}^k\sigma_{\curvearrowleft(i_s,j_s)}
\bigsqcup\coprod\limits_{s=1}^k S_{(i_s,j_s)\looparrowright}(\sigma)
\bigsqcup\coprod\limits_{s=1}^kS_{(i_s,j_s)\upharpoonleft\downharpoonright}(\sigma)$
\item[(ii)] Let $\sigma'\in \bS_n^2(k)$ be such that
$\sigma'\prec\sigma.$ Then $\sigma'\in \D(\sigma)$ iff
$\codim_{\ov\B_\sigma}\B_{\sigma'}=1.$ \end{itemize} \end{thm}

\subsection{}\label{3.8}
Now we are ready to finish the description of $\C(\sigma).$ Note
first that
$$\D(\sigma)\st\C(\sigma)\st \D(\sigma)\cup\N(\sigma)$$.

So we have to  find $ \N(\sigma)\cap \C(\sigma).$

Consider $\sigma=(i_1,j_1)\ldots(i_k,j_k)\in\bS_n^2.$ Let
$E_{\max}(\sigma)$ be the subset of arcs in $E(\sigma)$ without
fixed points outside of them. In that case all the fixed points of
$P_\sigma$ are under arc $(i,j)$ so that $(i,j)\in E_{\max}(\sigma)$
iff $f'_{(i,j)}(P_\sigma)=n-2k.$ We can use this fact to define
$E_{\max}(\sigma)$ formally:
$$E_{\max}(\sigma): =\{(i,j)\in E(\sigma)\ |\
f'_{(i,j)}(P_\sigma)=n-2k\}$$
\begin{lem} Consider $\sigma=(i_1,j_1)\ldots(i_k,j_k)\in \bS_n^2$
and let $(i,j)\in E(\sigma).$
\begin{itemize}
 \item[(i)] If $(i,j)\not\in E_{\max}(\sigma)$ then
 $\sigma_{(i,j)}^-\not\in\C(\sigma).$
 \item[(ii)] If $(i,j)\in E_{\max}(\sigma)$ then
$\dim \B_{\sigma_{(i,j)}^-}=\dim\B_\sigma-1.$
\end{itemize}
 \end{lem}
 \begin{proof}
  To prove (i)
 note that for any $(i,j)\in E(\sigma)$ one has ${\bf
over}_{(i,j)}(P_\sigma)=\emptyset.$ So the conditions of \ref{3.4}
for $(i,j)\in E(\sigma)$ are satisfied iff $(i,j)\not\in
E_{\max}(\sigma)$ and in this case
$S=\{\sigma_{\curvearrowleft(i,j)},\
\sigma_{(i,j)\curvearrowright}\}\ne\emptyset.$  Let $\sigma'\in S$
and let $P_{\sigma'}$ be obtained from $P_\sigma$ by replacing
$(i,j)$ to $(i',j')$. Then obviously
$\sigma_{(i,j)}^-=(\sigma')_{(i',j')}^-$ so that
$\sigma_{(i,j)}^-\prec\sigma'\prec\sigma$ and
$\sigma_{(i,j)}^-\not\in \C(\sigma).$

  To compute (ii) we use notation from the proof of theorem
\ref{3.1}. Consider an arc $(i_s,j_s)\in P_\sigma.$
\begin{itemize}
\item{} If $(i_s,j_s)$ does not intersect $(i,j)$ then
$$f'_{(i_s,j_s)}(P_\sigma)=f'_{(i_s,j_s)}(P_{\sigma_{(i,j)}^-})\quad{\rm
and}\quad
|c_{(i_s,j_s)}^l(P_\sigma)|=|c_{(i_s,j_s)}^l(P_{\sigma_{(i,j)}^-})|.
$$
\item{} If $(i_s,j_s)$ intersects $(i,j)$ on the left  $i$ is a new
fixed point under $(i_s,j_s)$ in $P_{\sigma_{(i,j)}^-}$ comparing to
$P_\sigma$ and the crossings from the left of $(i_s,j_s)$ are the
same in both patterns, so that
$$f'_{(i_s,j_s)}(P_\sigma)=f'_{(i_s,j_s)}(P_{\sigma_{(i,j)}^-})-1\quad{\rm
and}\quad
|c_{(i_s,j_s)}^l(P_\sigma)|=|c_{(i_s,j_s)}^l(P_{\sigma_{(i,j)}^-})|
$$
\item{} If $(i_s,j_s)$ intersects $(i,j)$ on the right then $j$ is a new
fixed point under $(i_s,j_s)$ in $P_{\sigma_{(i,j)}^-}$ comparing to
$P_\sigma$ and there is one less crossing from the left of
$(i_s,j_s)$ in $P_{\sigma_{(i,j)}^-}$ comparing to $P_\sigma$, so
that
$$f'_{(i_s,j_s)}(P_\sigma)=f'_{(i_s,j_s)}(P_{\sigma_{(i,j)}^-})-1\quad{\rm
and}\quad
|c_{(i_s,j_s)}^l(P_\sigma)|=|c_{(i_s,j_s)}^l(P_{\sigma_{(i,j)}^-})|+1
$$
\item{}
Finally, since there are no fixed points outside of $(i,j)$ in
$P_\sigma$ one has $f'_{(i,j)}(P_\sigma)=n-2k$
\end{itemize}
 Summarizing we get
$$\begin{array}{ll}
f(P_\sigma)&=\sum\limits_{s=1}^k
f'_{(i_s,j_s)}(P_\sigma)=\sum\limits_{(i_s,j_s)\ne(i,j)}
f'_{(i_s,j_s)}(P_\sigma)+f'_{(i,j)}(P_\sigma)=\\
&\sum\limits_{(i_s,j_s)\ne(i,j)}f'_{(i_s,j_s)}(P_{\sigma_{(i,j)}^-})-
|c_{(i,j)}^l(P_\sigma)|-|c_{(i,j)}^r(P_\sigma)|+n-2k=\\
&f(P_{\sigma_{(i,j)}^-})-
|c_{(i,j)}^l(P_\sigma)|-|c_{(i,j)}^r(P_\sigma)|+n-2k;\\
\end{array}
$$
$$
\begin{array}{l}
c(P_\sigma)=\sum\limits_{s=1}^k
|c_{(i_s,j_s)}^l(P_\sigma)|=\sum\limits_{(i_s,j_s)\ne(i,j)}
|c_{(i_s,j_s)}^l(P_\sigma)|+ |c_{(i,j)}^l(P_\sigma)|=\\
\left(\sum\limits_{{(i_s,j_s)\ne(i,j)}\atop{(i_s,j_s)\not\in
c_{(i,j)^r}(P_\sigma)}}|c_{(i_s,j_s)}^l(P_{\sigma_{(i,j)}^-})|\right)+
\left(\sum\limits_{(i_s,j_s)\in
c_{(i,j)^r}(P_\sigma)}|c_{(i_s,j_s)}^l(P_{\sigma_{(i,j)}^-})|+
|c_{(i,j)}^r(P_\sigma)|\right) +|c_{(i,j)}^l(P_\sigma)|=\\
\quad\\
 c(P_{\sigma_{(i,j)}^-})+|c_{(i,j)}^r(P_\sigma)| +|c_{(i,j)}^l(P_\sigma)|\\
\end{array}$$
Thus,
$$\begin{array}{ll}
\dim\B_\sigma&=k(n-k)-f(P_\sigma)-c(P_\sigma)=\\
 & k(n-k)-\left(
f(P_{\sigma_{(i,j)}^-})-
|c_{(i,j)}^l(P_\sigma)|-|c_{(i,j)}^r(P_\sigma)|+n-2k\right)-\\
&\quad\left( c(P_{\sigma_{(i,j)}^-})+|c_{(i,j)}^r(P_\sigma)|
+|c_{(i,j)}^l(P_\sigma)|\right)=\\
&k(n-k)-(n-2k)-f(P_{\sigma_{(i,j)}^-})-c(P_{\sigma_{(i,j)}^-})=\\
&(k-1)(n-k+1)+1-f(P_{\sigma_{(i,j)}^-})-c(P_{\sigma_{(i,j)}^-})=\dim\B_{\sigma_{(i,j)}^-}+1.\\
\end{array}
$$

\end{proof} \subsection{}\label{3.9} As a straightforward corollary
of theorem \ref{3.7} and lemma \ref{3.8}, we get \begin{thm} Let
$\sigma=(i_1,j_1)\ldots(i_k,j_k)\in\bS_n^2.$ One has \begin{itemize}
\item[(i)] $\C(\sigma)=\D(\sigma)\bigsqcup\{\sigma_{(i,j)}^-\ |\
(i,j)\in E_{\max}(\sigma)\}.$ \item[(ii)] $\sigma'\in \C(\sigma)$
iff $\codim_{\ov\B_\sigma}\B_{\sigma'}=1.$ \end{itemize} \end{thm}

\subsection{}\label{3.9a} In section 5 we will discuss intersections
of $\bB$ orbit closures in $\X_n.$ As a preparation we would like to
finish this section with a few preliminary general definition and
results.

 Given  $\emptyset\ne S\st\bS_n^2.$ We
define the components of $U=\bigcap\limits_{\sigma\in
S}\ov\B_\sigma$ in a standard way by setting
$\B_{\tau_1},\ldots,\B_{\tau_r}$ to be the components of $U$ if
$U=\bigcup\limits_{i=1}^r\ov\B_{\tau_i}$ and $\tau_i\preceq\tau_j$
implies $i=j$ for any $i,j\ :\ 1\leq i,j\leq r.$ Let us call
$\tau_1,\ldots,\tau_i$ component involutions of $U.$

Given $\sigma\in\bS_n^2(k)$ where $k\leq {1\over 2}n-1$. For  $m\ :\ k<m\leq {\frac 12}n$ put $F_m(\sigma):=\{\tau\in\bS_n^2(m)\ |\ \tau\succ\sigma\}.$ Let
$\{p_i\}_{i=1}^{n-2k}$ be fixed points of $\sigma$. Note that
$n-2k\geq 2(m-k)$ so that
$\sigma'=\sigma(p_1,p_2)\ldots(p_{2(m-k)-1},p_{2(m-k)}\in
\bS_n^2(m)$ and $\sigma'\succ\sigma.$ Thus,
$F_m(\sigma)\ne\emptyset$

As an easy corollary of subsection \ref{2.5} we get

\begin{prop} Let $k\leq {1\over 2}n-1$. For any
$\sigma\in\bS_n^2(k)$ and any $m\ :\ k<m\leq {1\over 2}n$ there
exists the unique $\ov\sigma_m\in F_m(\sigma)$ such that
$$\bigcap\limits_{\tau\in
F_m(\sigma)}\ov\B_\tau=\ov\B_{\ov\sigma_m}.$$

 In particular, let $\emptyset\ne S\st\bS_n^2(k)$ and let
$U=\bigcap\limits_{\sigma\in S}\ov\B_\sigma.$ Then all the component
involutions of $U$ are of length $k.$
  \end{prop}
   \begin{proof}
By $(*)$ of \ref{2.5} $U'=\bigcap\limits_{\tau\in
F_m(\sigma)}\ov\B_\tau\bigcap\Or_{(n-m,m)^*}\ne\emptyset.$ Let
$\tau_1,\ldots\tau_r$ be the component involutions of $U'.$  Then
by irreducibility of $\B_\sigma$ there exists $\tau_i=\ov\sigma_m$
such that $\B_\sigma\st\ov\B_{\ov\sigma_m}$ but then
$\ov\sigma_m\in F_m(\sigma)$ and $$U'= \bigcup\limits_{j=1}
^r\left(\ov\B_{\tau_j}\cap\ov\B_{\ov\sigma_m}\cap\Or_{(n-m,m)^*}\right)=\ov\B_{\ov\sigma_m}\cap\Or_{(n-m,m)^*}.$$
To complete the proof we have to show that
$\bigcap\limits_{\tau\in
F_m(\sigma)}\ov\B_\tau=\ov\B_{\ov\sigma_m}$ (in other words we
have to show that $\ov{\bigcap\limits_{\tau\in
F_m(\sigma)}\ov\B_\tau\bigcap\Or_{(n-m,m)^*}}=\bigcap\limits_{\tau\in
F_m(\sigma)}\ov\B_\tau\bigcap\ov\Or_{(n-m,m)^*}=\bigcap\limits_{\tau\in
F_m(\sigma)}\ov\B_\tau$, where the second equality is obvious). To do this we must show that for any
$j<m$ for any $\B_{\sigma'} \st \bigcap\limits_{\tau\in
F_m(\sigma)}\ov\B_\tau\bigcap\Or_{(n-j,j)^*}$ one has
$\B_{\sigma'}\st\ov\B_{\ov\sigma_m}.$  Let us  consider
$F_m(\sigma').$ Since $\B_{\sigma'}\st\bigcap\limits_{\tau\in
F_m(\sigma)}\ov\B_\tau$   we get $F_m(\sigma)\st F_m(\sigma')$. By
the previous there exists the unique minimal element
$\ov\sigma'_m$ of $F_m(\sigma')$ so that
$\ov\sigma'_m\preceq\ov\sigma_m$ which provides
$\B_{\sigma'}\st\ov\B_{\ov\sigma_m}$ and completes the proof.

Let us consider $U$ from the proposition. It is enough to show that
for any $\sigma$ such that $\B_\sigma\st U$ and $l(\sigma)<k$ there
exists $\sigma'\succ\sigma$ such that $\B_{\sigma'}\st U$ and
$l(\sigma')=k.$ Consider $F_k(\sigma).$ By the definition of $U$ we
get $S\st F_k(\sigma)$ so that  $\ov\sigma_k$  satisfies
 $\B_{\ov\sigma_k}\st U.$
\end{proof}

 Note that in general $U$ defined in the proposition can be
 reducible.

\subsection{}\label{3.10} Given $\sigma'\in F_p(\sigma)$, to show
that $\sigma'=\ov\sigma_p$ it is enough to check that
$$\dim\B_{\sigma^{\prime\prime}}\geq \dim\B_{\sigma'}\quad {\rm for\
any}\quad \sigma^{\prime\prime}\in F_p(\sigma).\eqno{(*)}$$

 It is easy to write explicitly $\ov\sigma_{k+1}$ and to compute the dimension of
 $\B_{\ov\sigma_{k+1}}.$
 For $\sigma\in\bS_n^2(k)$ let
$i_\sigma:=\min\{\{s\}_{s=1}^n\setminus\langle\sigma\rangle\}$ and
$j_\sigma:=\max\{\{s\}_{s=1}^n\setminus\langle\sigma\rangle\}$. If
$E_{\max}(\sigma)\ne\emptyset$ then $(i,j)\in{\rm
over}_{i_\sigma}(P_\sigma),{\rm over}_{j_\sigma}(P_\sigma)$ for any
$(i,j)\in E_{\max}$. Let $E_{\max}(\sigma)=\{(i_r,j_r)\}_{r=1}^s$
where the arcs are ordered in the increasing order according to the
left end points. Note that
$i_1<i_2<\ldots<i_s<i_\sigma<j_\sigma<j_1<\ldots j_s.$   Put
$j_0:=j_\sigma$ and $i_{s+1}:=i_\sigma$ and let $\sigma'$ be
obtained from $\sigma$ by deleting all the cycles in
$E_{\max}(\sigma).$ Put
$$\overline\sigma_{k+1}=\left\{\begin{array}{ll}
\sigma(i_\sigma,j_\sigma)& {\rm if}\ E_{\max}(\sigma)=\emptyset\\
\sigma'(i_1,j_0)(i_2,j_1)\ldots(i_s,j_{s-1})(i_{s+1},j_s)& {\rm otherwise}\\
\end{array}\right.$$
Let us show some examples to make the picture clear:
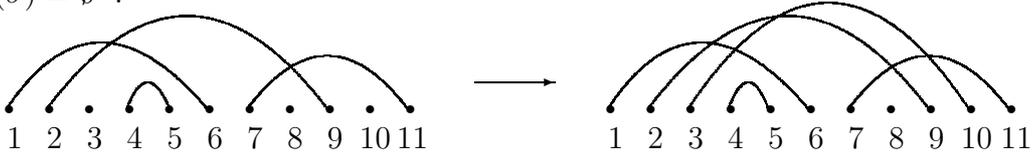
\begin{center}
\begin{picture}(100,100)\put(-170,80){$E_{\max}(\sigma)=\emptyset\ :\ $}
\multiput(-139,40)(15,0){11}%
{\circle*{3}}
 \put(-140,25){1}
 \put(-125,25){2}
 \put(-110,25){3}
 \put(-95,25){4}
 \put(-80,25){5}
 \put(-65,25){6}
 \put(-50,25){7}
 \put(-35,25){8}
\put(-20,25){9}
 \put(-8,25){10}
 \put(6,25){11}

 \qbezier(-140,40)(-107,90)(-65,40)
 \qbezier(-125,40)(-72.5,110)(-20,40)
\qbezier(-50,40)(-20,80)(10,40)
 \qbezier(-95,40)(-87.5,60)(-79,40)

%\put(35,60){$\sigma(3,10)$}
 \put(35,50){\vector(1,0){30}}

\multiput(86,40)(15,0){11}%
{\circle*{3}}
 \put(85,25){1}
 \put(100,25){2}
 \put(115,25){3}
 \put(130,25){4}
 \put(145,25){5}
 \put(160,25){6}
 \put(175,25){7}
 \put(190,25){8}
\put(205,25){9}
 \put(217,25){10}
 \put(232,25){11}

 \qbezier(85,40)(117.5,90)(160,40)
 \qbezier(100,40)(152.5,110)(205,40)
\qbezier(175,40)(205,80)(235,40)
 \qbezier(130,40)(137.5,60)(145,40)
\qbezier(115,40)(167.5,120)(220,40)
\end{picture}
\end{center}
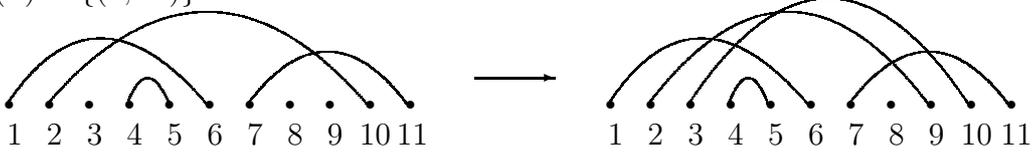
\begin{center}
\begin{picture}(100,100)\put(-170,80){$E_{\max}(\sigma)=\{(2,10)\}\ :\ $}
\multiput(-139,40)(15,0){11}%
{\circle*{3}}
 \put(-140,25){1}
 \put(-125,25){2}
 \put(-110,25){3}
 \put(-95,25){4}
 \put(-80,25){5}
 \put(-65,25){6}
 \put(-50,25){7}
 \put(-35,25){8}
\put(-20,25){9}
 \put(-8,25){10}
 \put(6,25){11}

 \qbezier(-140,40)(-107,90)(-65,40)
 \qbezier(-125,40)(-65,110)(-5,40)
\qbezier(-50,40)(-20,80)(10,40)
 \qbezier(-95,40)(-87.5,60)(-79,40)

%\put(35,60){$\sigma(3,10)$}
 \put(35,50){\vector(1,0){30}}

\multiput(86,40)(15,0){11}%
{\circle*{3}}
 \put(85,25){1}
 \put(100,25){2}
 \put(115,25){3}
 \put(130,25){4}
 \put(145,25){5}
 \put(160,25){6}
 \put(175,25){7}
 \put(190,25){8}
\put(205,25){9}
 \put(217,25){10}
 \put(232,25){11}

 \qbezier(85,40)(117.5,90)(160,40)
 \qbezier(100,40)(152.5,110)(205,40)
\qbezier(175,40)(205,80)(235,40)
 \qbezier(130,40)(137.5,60)(145,40)
\qbezier(115,40)(167.5,120)(220,40)
\end{picture}
\end{center}

 Note that
$$(R_{\ov\sigma_{k+1}})_{i,j}=\left\{\begin{array}{ll}
(R_\sigma)_{i,j}+1&{\rm if}\ i\leq i_r\ {\rm and}\ j\geq j_{r-1}\
{\rm
for}\ 1\leq r\leq s\\
(R_\sigma)_{i,j}& {\rm otherwise}\\
\end{array}\right.$$
So that $\ov\sigma_{k+1}\succ\sigma.$ We must show that
$\ov\sigma_{k+1}$ defined above is indeed the minimal element of
$F_{k+1}(\sigma).$
\begin{prop}
 Let $k\leq {1\over 2}n-1$.  Then for every $\sigma'\in F_p(\sigma)$ where $p\geq k+1$
one has $\sigma'\succeq \ov\sigma_{k+1}.$ One also has
$\codim_{\ov\B_{\ov\sigma_{k+1}}}\B_\sigma=|E_{\max}(\sigma)|+1.$
\end{prop}
\begin{proof}
Let us show this by induction on the $|E_{\max}(\sigma)|$. Put
$A(\sigma):=\{\sigma'\ |\ \sigma\in \D(\sigma')\}$ and
$B(\sigma):=\{\sigma'\ |\ \sigma\in \C(\sigma')\ {\rm and}\
\ell(P_{\sigma'})=\ell(P_\sigma)+1\}.$ Note that for any $\sigma'\in
A(\sigma)$ one has
$$|E_{\max}(\sigma)|-1\leq |E_{\max}(\sigma')|\leq
|E_{\max}(\sigma)|.\eqno(**)$$

Note also that by definition
$$B(\sigma)=\left\{\begin{array}{ll}\sigma(i_\sigma,j_\sigma)& {\rm if}\
E_{\max}(\sigma)=\emptyset\\
\emptyset&{\rm otherwise}\\
\end{array}\right.\eqno{({**\atop *})}$$

If $E_{\max}(\sigma)=\emptyset$ then
$B(\sigma)=\{\sigma(i_\sigma,j_\sigma)\}$ so that
$\codim_{\ov\B_{\sigma(i_\sigma,j_\sigma)}}\B_\sigma=1.$ For any
$\sigma'\in F_{k+1}(\sigma)$ one has
$\codim_{\ov\B_{\sigma'}}\B_\sigma\geq 1$, thus, by $(*)$  and
Proposition \ref{3.10} we get the result.

Now assume this is true for $\sigma\in\bS_n^2(k)$ such that
$|E_{\max}(\sigma)|=s$  and show this is true for $\sigma$ such that
$|E_{\max}(\sigma)|=s+1$. Let $\sigma'$ be obtained from $\sigma$ by
changing $(i_1,j_1)$ to $(i_1, j_\sigma).$ Note that
$\sigma=\sigma'_{(i_1,j_\sigma)\curvearrowright}$ so that
$\sigma'\in A(\sigma)$ and that $|E_{\max}(\sigma')|=s.$ Note also
that $\ov\sigma_{k+1}=\ov\sigma'_{k+1}$, thus, by induction one has
$\codim_{\ov\B_{\ov\sigma_{k+1}}}\B_\sigma=\codim_{\ov\B_{\ov\sigma_{k+1}}}\B_{\sigma'}
+\codim_{\ov\B_{\sigma'}}\B_\sigma=s+1.$

To show that $\ov\sigma_{k+1}$ is the minimal element of
$F_{k+1}(\sigma)$ it is enough to show by $(*)$ that
$\codim_{\ov\B_{\sigma'}}\B_\sigma\geq s+1$ for any $\sigma'\in
F_{k+1}(\sigma).$ Assume this is not true. Let $\sigma\in\bS_n^2(k)$
be maximal such that $|E_{\max}(\sigma)|=s+1$ and there exists
$\sigma'\in F_{k+1}(\sigma) $ such that
$\codim_{\ov\B_{\sigma'}}\B_\sigma\leq s.$ By $({**\atop *})$
$B(\sigma)=\emptyset$ so that by theorem \ref{3.9} there exists
$\sigma^{\prime\prime}\in A(\sigma)$ such that
$\sigma'\succ\sigma^{\prime\prime}\succ\sigma.$ Note that
$\sigma'\in F_{k+1}(\sigma^{\prime\prime}).$ By $(**)$ $s\leq
|E_{\max}(\sigma^{\prime\prime})|\leq s+1$. If
$|E_{\max}(\sigma^{\prime\prime})|=s+1$ we get that
$\codim_{\ov\B_{\sigma'}}\B_{\sigma^{\prime\prime}}\leq s-1$ which
contradicts maximality of $\sigma.$ If
$|E_{\max}(\sigma^{\prime\prime})|=s$ then by induction hypothesis
$\codim_{\ov\B_{\sigma'}}\B_{\sigma^{\prime\prime}}\geq s$ so that
$\codim_{\ov\B_{\sigma'}}\B_{\sigma}\geq s+1$ which contradicts the
assumption.
\end{proof}

\section{Combinatorial description of the closure of an
orbital variety of nilpotent order 2}

\subsection{}\label{4.1}
In this section we apply the results obtained above to the closures
of $\bB$ orbits of dimension $k(n-k)$ in $\Bb_{(n-k,k)}$, that is of
the maximal possible dimension.  Such orbit is dense in an orbital
variety associated to $\Or_{(n-k,k)^*}$ and every orbital variety
associated to $\Or_{(n-k,k)^*}$ admits such $\bB$ orbit as we
explained in \ref{2.2}.

 As it was shown in \cite[4.5]{Mo2} the
closure of an orbital variety of nilpotent order 2 is a union of
orbital varieties. In \cite[4.3]{Mo2} the combinatorial description
of such orbital variety in terms of Young tableaux is given. We
would like to use link patterns to simplify this description.

By Proposition \ref{2.2}   $\B_T$ is the dense $\bB$-orbit in $\V_T$
for $T\in\bT_{(n-k,k)^*}.$ By definition of an orbital variety
$\ov\V_T=\V_T\bigsqcup\bigcup\limits_{\sigma'\in\N(\sigma_T)}\ov\B_{\sigma'}$.
Moreover, since an element of $\N(\sigma_T)$ is obtained from
$\sigma_T$ by erasing an external arc it is a link pattern without
fixed points under an arc and without arc crossings, thus,
$\B_{\sigma'}$ for $\sigma'\in \N_{\sigma_T}$ is
$\sigma'=\sigma_{T'}$ for some $T'\in\bT_{(n-k+1,k-1)^*}.$

How to define $\{T'\in \bT_{(n-k+1,k-1)^*}\ |\
\sigma_{T'}\in\N(\sigma_T)\}$? The answer given in \cite[4.3]{Mo2}
is too complex to be satisfactory. Let us try to understand the
picture using link patterns.

Given  $\sigma_T=(i_1,j_1)\ldots (i_k,j_k)$  ordered in increasing
order according to the second entry of the cycle (that is
$j_1<\ldots<j_k$). Then $T_2=\left(\begin{array}{c}j_1\cr\vdots\cr
j_k\cr\end{array}\right)$ (cf. \ref{2.2}). Moreover, set
$\{j_1,\ldots j_k\}$ defines $\sigma_T$ completely. In turn
$\sigma_{T'}$ is obtained from $\sigma_T$ by erasing some external
arc, that is $T'_2=\left(\begin{array}{c}\vdots\\j_{l-1}\\j_{l+1}\\
\vdots \cr\end{array}\right)$  for some $l\ :\ 1\leq l\leq k.$ To
determine $T'\ :\ \sigma_{T'}\in\N_{\sigma_T}$ we must find
$E(\sigma_T)$ -- the subset of external arcs.

Note that since arcs of $P_{\sigma_T}$ do not intersect, for any $t>s$ (that is such that $j_t>j_s$)
one has $i_t\not\in[i_s,j_s]$ that is either $i_t<i_s$ or $i_t>j_s.$  Thus,  arc $(i_s,j_s)\in
E(\sigma_T)$  if for any $(i_t,j_t)\in\sigma$ such that $t>s$ one
has $i_t>j_s$. Hence, for
any $t>s$ there must be at least $2(t-s)-1$ points between $j_s$ and
$j_t.$ In other words the necessary condition is $j_t-j_s\geq
2(t-s).$ This is also a sufficient condition. Let us formulate this statement as
a lemma
\begin{lem}
Let $T\in\bT_{(n-k,k)^*}$ be such that
$T_2=\left(\begin{array}{c}j_1\cr\vdots\cr
j_k\cr\end{array}\right).$ Let $\sigma_T=(i_1,j_1)\ldots(i_k,j_k).$
Then $E(\sigma_T)=\{(i_l,j_l)\ |\ l=k\ or\ \forall\ s\ :\ l<s\leq k\
\Rightarrow\ j_s-j_l\geq 2(s-l)\}.$
\end{lem}
\begin{proof}
First of all let us show that if $(i_l,j_l)\in E(\sigma_T)$ then for
any $s>l$ one must have $j_s-j_l\geq 2(s-l).$  Indeed, for any $t\
:\ l<t\leq s$ one has (as we already mentioned) $i_t>j_l$ and
$j_t\leq j_s$ so that $[j_l+1,j_s]$ must contain at least $2(s-l)$
points, therefore $j_s-j_l\geq 2(s-l).$

Exactly in the same way we show that if  $j_s-j_l\geq 2(s-l)$ for
any $s>l$  then $(i_l,j_l)\in E(\sigma_T).$ Indeed, for $s=l+1$ one
has $j_{l+1}\geq j_l+2$ so that $j_{l+1}-1\in\langle T_1\rangle$ and
$i_{l+1}=j_{l+1}-1>j_l.$ Assume that we already know arcs
$(i_{l+1},j_{l+1}),\ldots (i_{s-1},j_{s-1})$  and
$i_{l+1},\ldots,i_{s-1}\in [j_l+1,j_{s-1}].$ By the construction
$i_s=\max\{r\}_{r=1}^n\setminus \{i_p,j_p\}_{p=1}^{s-1}$ And since
there are at least $2(s-l)$ points in $[j_l+1,j_s]$ and only
$2(s-1-l)+1$ are occupied we get that $i_s\in[j_l+1,j_s].$
\end{proof}

\subsection{}\label{4.2} Consider $T\in\bT_{(n-k,k)^*}$ and let
$T_2=\left(\begin{array}{c}j_1\cr\vdots\cr
j_k\cr\end{array}\right).$ For $b\in \langle T_2\rangle$ put
$T\langle b\rangle$ to be the tableau obtained from $T$ by moving
$b$ from $T_2$ to an appropriate place in $T_1.$

Let us denote by $\N(T):=\{S\in\bT_{(n-k+1,k-1)^*}\ |\ \V_S\st
\V_T\}$. By dimension reasoning it is obvious that $\V_S\ :\
S\in\N(T)$ are orbital varieties of maximal dimension in
$\ov\V_T\setminus \V_T.$ Moreover the fact that $\V_T$ admits a
dense $\bB$ orbit and Proposition \ref{3.3} provide that these are
all maximal orbital varieties in $\ov\V_T\setminus \V_T.$ In other
words $\ov\V_T=\V_T\bigsqcup\bigcup\limits_{S\in\N(T)}\ov\V_S.$

 As a straightforward corollary of Lemma \ref{4.1} we get
\begin{thm}
Let $T\in\bT_{(n-k,k)^*}$ and
$T_2=\left(\begin{array}{c}j_1\cr\vdots\cr
j_k\cr\end{array}\right)$. Then
$$\N(T)=\{T\langle j_i\rangle\ |\ i=k\ or\ \forall\ s\ :\
i<s\leq k\ \Rightarrow\ j_s-j_i\geq 2(s-i)\}.$$
\end{thm}
This theorem simplifies Theorem 4.3 from \cite{Mo2}.

\section{Meanders and intersections of the components of a
Springer fiber of nilpotent order 2}
\subsection{}\label{4.3}
In this section we apply our results to the theory of intersections
of the components of Springer fiber $\F_x$ where $x\in
\Or_{(n-k,k)^*}.$ For $S,T\in\bT_{(n-k,k)^*}$ let $\F_T,\ \F_S$ be
the corresponding components of $\F_x.$ As we explained in \ref{1.8}
the number  of components of $\F_T\cap\F_S$ and their codimensions
equal to the number of components of $\V_T\cap\V_S$ and their
codimensions. So we formulate our results in terms of the components
of a Springer fiber, however the computations are made in $\bB$
orbits of nilpotent order 2. First note that for
$T,S\in\bT_{(n-k,k)^*}$ one has
$$\V_T\cap\V_S=\ov\B_T\cap\ov\B_S\cap\Or_{(n-k,k)^*}=
\coprod\limits_{{\sigma\in\bS_n^2(k)}\atop {\sigma\preceq
\sigma_T,\sigma_S}}\B_\sigma.$$

 We will try to understand these intersections in terms of meanders. We draw one link pattern
 upwards and one link pattern downwards. If both link patterns of a meander are without
 intersections we get
 a classical meander, otherwise we get a generalized meander.

  As we have
mentioned in \ref{1.8}  a generalized  meander $M_{\sigma,\sigma'}$
corresponding to $\sigma,\sigma'\in\bS_n^2(k)$ provides us an easy
way to compute $R_{\sigma,\sigma'}$. However for a generalized
meander the picture seems to be too complex to understand it without
going to  matrix $R_{\sigma,\sigma'}$.

As it is shown in \cite[3.8]{M-PII} even an intersection of
codimension 1 for (closures of) generic $\bB$ orbits of the same dimension can be
reducible. Moreover, it can be not of pure dimension as it is shown
by the example below. We also use this example to construct the corresponding
generalized meander. Consider $\sigma=(1,3)(4,5),
\ \sigma'=(2,3)(4,6)\in\bS_6^2(2).$ Note that
$\dim\B_{(1,3)(4,5)}=\dim\B_{(2,3)(4,6)}=2\cdot 4-1=7.$ The generalized meander
is
\begin{center}
\begin{picture}(100,80)
\put(-70,40){$M_{(1,3)(4,5),\ (2,3)(4,6)}=$}
\multiput(40,40)(20,0){6}%
{\circle*{3}}

 \qbezier(40,40)(60,70)(80,40)
 \qbezier(100,40)(110,60)(120,40)
\qbezier(60,40)(70,20)(80,40)
 \qbezier(100,40)(120,10)(140,40)
\end{picture}
\end{center}
and
$$R_{(1,3)(4,5),\ (2,3)(4,6)}=\left(\begin{array}{cccccc}
0&0&1&1&1&2\\
0&0&0&0&1&1\\
0&0&0&0&0&1\\
0&0&0&0&0&1\\
0&0&0&0&0&0\\
0&0&0&0&0&0\\
\end{array}\right)$$
One can easily check using  Proposition \ref{2.3} that
$R_{(1,3)(4,5),\ (2,3)(4,6)}\not\in\bR_6^2$ so that the intersection
is reducible by Theorem \ref{2.5}.  Computing $R_{(1,3)(4,6)}$ and
$R_{(1,6)(2,5)}$ we get
$\overline\B_{(1,3)(4,5)}\cap\overline\B_{(2,3)(4,6)}=\overline\B_{(1,3)(4,6)}\cup
\overline\B_{(1,6)(2,5)}$ so that the intersection is reducible and $\dim\B_{(1,3)(4,6)}=2\cdot 4-2=6$, $\dim\B_{(1,6)(2,5)}=2\cdot 4-4=4.$

\subsection{}\label{4.4}
 As we explained in the introduction the
meanders consisting of link patterns of orbital varieties that is
$M_{\sigma_S,\sigma_T}$, where $S,T\in\bT_{(n-k,k)^*}$ are the most
simple on one hand and the most important on the other hand.
To simplify the notation
we put $P_T:=P_{\sigma_T}$ and  $M_{S,T}:=M_{\sigma_S,\sigma_T}.$

By the very nature of the ordering of $\bB$ orbits if $\V_T\cap\V_S$
is of higher (than 1) codimension and $\B_\sigma$ is a component of the intersection
there exist $\B_\tau\st\V_T,\ \B_\tau\not\st\V_S$
and $\B_{\tau'}\st\V_S,\ \B_{\tau'}\not\st\V_T$ such that $\dim\B_\tau=\dim\B_{\tau'}=\dim\B_\sigma+1$ and $\B_\sigma$ is a component of their
intersection. Thus, the description of the intersections of higher codimension  involves the description of intersections of codimension 1 for general $\bB$ orbits.
Just to illustrate this point let us show that the example of \ref{4.3} provides us with the intersection of codimension 2 of two orbital  varieties which is not only reducible but also is not of pure dimension. This is the example from  \cite[5.7]{M-PI}. Consider
$$T=\begin{array}{ll}
1&3\\
2&6\\
4& \\
5& \\
\end{array}\quad{\rm and}\quad S=\begin{array}{ll}
1&2\\
3&5\\
4& \\
6& \\
\end{array}$$
so that
$\sigma_T=(2,3)(5,6)$ and $\sigma_S=(1,2)(4,5).$ One has
$\dim\B_{(2,3)(5,6)}=\dim\B_{(1,2)(4,5)}=8.$ The corresponding
meander is
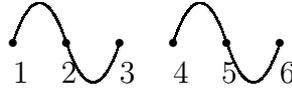
\begin{center}
\begin{picture}(100,80)
\multiput(-40,40)(20,0){6} {\circle*{3}}
 \put(-40,25){1}
 \put(-22,25){2}
 \put(0,25){3}
 \put(20,25){4}
 \put(38,25){5}
 \put(60,25){6}

 \qbezier(-40,40)(-30,70)(-20,40)
 \qbezier(20,40)(30,70)(40,40)
\qbezier(-20,40)(-10,10)(0,40)
 \qbezier(40,40)(50,10)(60,40)
 \end{picture}
 \end{center}
 Comparing with the meander in \ref{4.3} we see that $\B_{(1,3)(4,5)}\st\ov\B_{(1,2)(4,5)}$ of codimension 1 and $\B_{(1,3)(4,5)}\not\st\ov\B_{(2,3)(5,6)}$, also $\B_{(2,3)(4,6)}\st\ov\B_{(2,3)(5,6)}$ of codimension 1 and $\B_{(2,3)(4,6)}\not\st\ov\B_{(1,2)(4,5)}$
Moreover computing $R_{\sigma_T,\sigma_S}$ we get
 $$R_{\sigma_T,\sigma_S}=\left(\begin{array}{cccccc}
0&0&1&1&1&2\\
0&0&0&0&1&1\\
0&0&0&0&0&1\\
0&0&0&0&0&1\\
0&0&0&0&0&0\\
0&0&0&0&0&0\\\end{array}\right)=R_{(1,3)(4,5),\ (2,3)(4,6)}$$
Thus, exactly as in \ref{4.3} we get
 $$\V_T\cap\V_S=\left(\ov\B_{(1,3)(4,5)}\cap\ov\B_{(2,3)(4,6)}\right)\cap\Or_{(4,2)^*}=
 \left(\ov\B_{(1,3)(4,6)}\cup\ov\B_{(1,6)(2,5)}\right)\cap\Or_{(4,2)^*}$$
 and  $\dim\B_{(1,3)(4,6)}=6$ and $\dim\B_{(1,6)(2,5)}=4$
 so that $\V_T\cap\V_{S}$ contains one component of codimension 2
 and another component of codimension 4.

\subsection{}\label{4.5}
To formulate the statement on the intersections of codimension 1 we need more notation.

 Arcs of $M_{S,T}$ never intersect one another. A connected subset of arcs in
 $M_{S,T}$  can be either
 open or closed . If it is open we
 call it an  {\bf interval}. If it is closed we call it
 a {\bf
loop}. The number of arcs in a path $\Pa$ is called the length of
$\Pa$, we denote it by $L(\Pa)$.
Note that if $\Pa$ is a loop $L(\Pa)$ is even. An interval can
be either of even length, then we call it even, or of odd length, correspondingly
we call it odd.  We call a meander {\bf even} if it consists only of loops and even intervals. Otherwise we call it {\bf odd}. For example, the meander in \ref{4.4}
consists of two even intervals (each one of length 2) so it is an even meander.

Now we can formulate our result.
\begin{thm}
For $S,T\in\bT_{(n-k,k)^*}$ one has
 $\codim_{\F_T}(\F_T\cap\F_S)=1$ if and only if $M_{S,T}$ is an
 even  meander with $k-1$ loops. In this case $\F_T\cap\F_S$ is irreducible.
  \end{thm}
\begin{proof}
 Note
that $\codim_{\V_T}(\V_T\cap\V_S)=1$ iff
$\D(\sigma_T)\cap\D(\sigma_S)\ne\emptyset.$

Since for $\sigma'\in D(\sigma_T)$ one has
$\dim\B_{\sigma'}=k(n-k)-1$ we get that either there is one fixed
point under one arc of $P_{\sigma'}$  or exactly two arcs of
$P_{\sigma'}$ intersect.

Let $P_{\sigma'}$ be a link pattern with one fixed point under one
arc. Let $m$ be this fixed point and $(i,j)$ be the corresponding
arc. There are exactly two ways to get a bigger link pattern (that
is without fixed points under an arc and without crossings) and by
maximality of dimension each one of them is a link pattern of some
orbital variety:
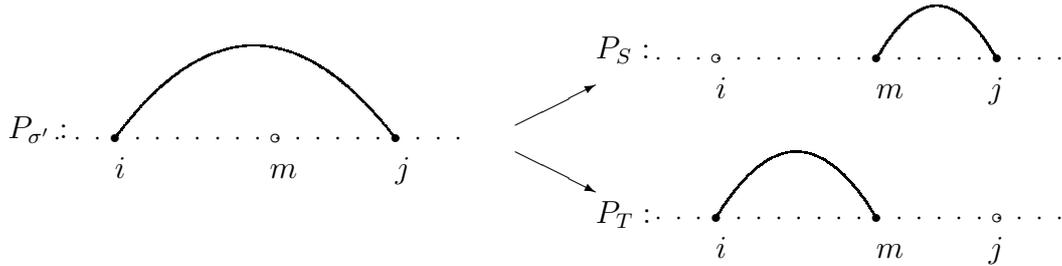
\begin{center}
\begin{picture}(100,140)
\put(-150,60){$P_{\sigma'}:$}
\multiput(-131.5,60)(7.5,0){21}%
{\circle*{1}}
 \put(-110,60){\circle*{3}}
 \put(-5,60){\circle*{3}}
 \put(-50,60){\circle{3}}
 \put(-110,45){$i$}
 \put(-52,45){$m$}
\put(-5,45){$j$}

 \qbezier(-110,60)(-58.5,130)(-5,60)

 \put(40,65){\vector(2,1){30}}
\put(70,90){$P_S:$}

\multiput(93.5,90)(7.5,0){21}%
{\circle*{1}}
 \put(115,90){\circle{3}}
 \put(175,90){\circle*{3}}
 \put(220,90){\circle*{3}}
 \put(115,75){$i$}
 \put(175,75){$m$}
 \put(217,75){$j$}

 \qbezier(175,90)(197.5,130)(220,90)

 \put(40,55){\vector(2,-1){30}}
 \put(70,28){$P_T:$}

 \multiput(93.5,30)(7.5,0){21}%
{\circle*{1}}
 \put(115,30){\circle*{3}}
 \put(175,30){\circle*{3}}
 \put(220,30){\circle{3}}
 \put(115,15){$i$}
 \put(175,15){$m$}
 \put(217,15){$j$}
 \qbezier(115,30)(145,80)(175,30)

\end{picture}
\end{center}
where dots on $[i+1,m-1]$ are all end points of arcs belonging to
$\pi_{i,m}(P_{\sigma'})$ and dots on $[m+1,j-1]$ are end points of
arcs belonging to $\pi_{m,j}(P_{\sigma'}).$

Note that $M_{S,T}$ in that case consists of $k-1$ loops
$(i_s,j_s)\ne(i,m),\ (m,j)$ (all the loops are of length 2) and a
unique even interval (of length 2) connecting $i$ and $j$ via $m$:
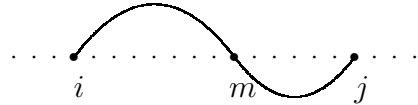
\begin{center}
\begin{picture}(100,40)
\multiput(-32.5,20)(7.5,0){21}%
{\circle*{1}}
 \put(-10,20){\circle*{3}}
 \put(95,20){\circle*{3}}
 \put(50,20){\circle*{3}}
 \put(-10,5){$i$}
 \put(95,5){$j$}
\put(48,5){$m$}

 \qbezier(-10,20)(20,60)(50,20)
 \qbezier(50,20)(72.5,-10)(95,20)
\end{picture}
\end{center}

Now let $P_{\sigma'}$ have two intersecting arcs. Let  $(i,j),
(i',j')$ where $i<i'<j<j'$ be these arcs. There are exactly two
ways to get a bigger link pattern (that is without crossings and
without fixed points under an arc) and by maximality of dimension
each one of them is a link pattern of some orbital variety:
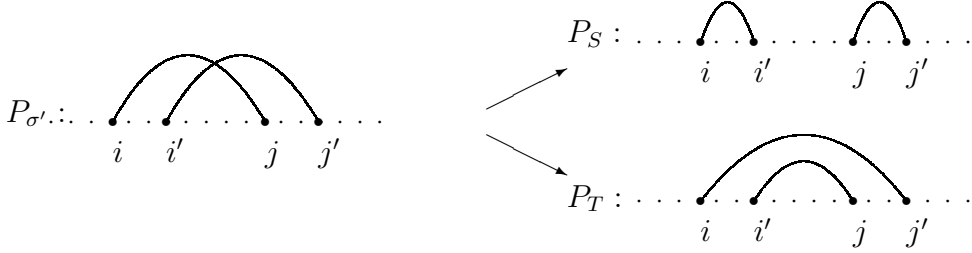
\begin{center}
\begin{picture}(100,100)
\put(-150,60){$P_{\sigma'}:$}
\multiput(-133,60)(7.5,0){3}%
{\circle*{1}}
 \put(-110,60){\circle*{3}}
 \put(-110,45){$i$}
 \multiput(-104,60)(7.5,0){3}%
{\circle*{1}}
 \put(-90,60){\circle*{3}}
 \put(-90,45){$i'$}
\multiput(-82.5,60)(7.5,0){4}%
{\circle*{1}}
 \put(-53,60){\circle*{3}}
 \put(-53,45){$j$}
 \put(-33,60){\circle*{3}}
 \multiput(-47,60)(7.5,0){3}%
{\circle*{1}}
 \put(-33,45){$j'$}
\multiput(-25,60)(7.5,0){3}%
{\circle*{1}}
 \qbezier(-110,60)(-82.5,110)(-53,60)
\qbezier(-90,60)(-62.5,110)(-33,60)

 \put(30,65){\vector(2,1){30}}
 \put(60,90){$P_S:$}

\multiput(87,90)(7.5,0){3}%
{\circle*{1}}
 \put(110,90){\circle*{3}}
 \put(110,75){$i$}
 \put(130,90){\circle*{3}}
 \put(130,75){$i'$}
\multiput(116,90)(7.5,0){2}%
{\circle*{1}}
\multiput(137.5,90)(7.5,0){4}%
{\circle*{1}}
 \put(167,90){\circle*{3}}
 \put(167,75){$j$}
 \put(187,90){\circle*{3}}
 \put(187,75){$j'$}
 \multiput(173,90)(7.5,0){2}%
{\circle*{1}}
\multiput(195,90)(7.5,0){3}%
{\circle*{1}}
 \qbezier(110,90)(120,120)(130,90)
\qbezier(167,90)(177,120)(187,90)

 \put(30,55){\vector(2,-1){30}}
\put(60,28){$P_T:$}
\multiput(87,30)(7.5,0){3}%
{\circle*{1}}
 \put(110,30){\circle*{3}}
 \put(110,15){$i$}
 \put(130,30){\circle*{3}}
 \put(130,15){$i'$}
 \multiput(116,30)(7.5,0){2}%
{\circle*{1}}
\multiput(137.5,30)(7.5,0){4}%
{\circle*{1}}
 \put(167,30){\circle*{3}}
 \put(167,15){$j$}
 \put(187,30){\circle*{3}}
 \put(187,15){$j'$}
 \multiput(173,30)(7.5,0){2}%
{\circle*{1}}
\multiput(195,30)(7.5,0){3}%
{\circle*{1}}
 \qbezier(110,30)(148.5,80)(187,30)
\qbezier(130,30)(148.5,60)(167,30)
\end{picture}
\end{center}
where dots on $[i+1,i'-1]$ are all end points of arcs belonging to
$\pi_{i,i'}(P_{\sigma'})$, dots on $[i'+1,j-1]$ are all end points
of arcs belonging to $\pi_{i',j}(P_{\sigma'})$ and dots on
$[j+1,j'-1]$ are all end points of arcs belonging to
$\pi_{j,j'}(P_{\sigma'}).$

Note that $M_{S,T}$ consists of $k-1$ loops: $k-2$ loops of type
$(i_s,j_s)$ (of length 2) where
$\{i_s,j_s\}\cap\{i,i',j,j'\}=\emptyset$ and a loop (of length 4)
involving $\{i,i',j,j'\}$:
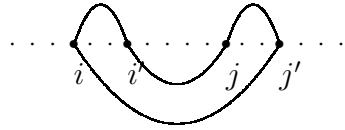
\begin{center}
\begin{picture}(100,80)
\multiput(-33,40)(7.5,0){3}%
{\circle*{1}}
 \put(-10,40){\circle*{3}}
 \put(-10,25){$i$}
 \put(10,40){\circle*{3}}
 \put(10,25){$i'$}
\multiput(-4,40)(7.5,0){2}%
{\circle*{1}}
\multiput(17.5,40)(7.5,0){4}%
{\circle*{1}}
 \put(47,40){\circle*{3}}
 \put(47,25){$j$}
 \put(67,40){\circle*{3}}
 \put(67,25){$j'$}
 \multiput(53,40)(7.5,0){2}%
{\circle*{1}}
\multiput(75,40)(7.5,0){3}%
{\circle*{1}}
 \qbezier(-10,40)(0,70)(10,40)
 \qbezier(47,40)(57,70)(67,40)
\qbezier(-10,40)(28.5,-20)(67,40)
 \qbezier(10,40)(28.5,10)(47,40)
 \end{picture}
\end{center}

Thus, in both cases we get an even meander with $k-1$ loops. Let
$M_{S,T}$ be a meander with $2k$ arcs. If it is even and contains
$k-1$ loops then it must contain $k-2$ loops of length 2 and
either one loop of length 4 or one more loop of length 2 and one
even interval of length 2, so that the cases above are all possible
cases for even meanders with $k-1$ loops.

 Also, one can see from the pictures above that for
$\sigma\in\D(\sigma_T)$ there exists exactly one $S$ such that
$\D_{\sigma_T}\cap\D_{\sigma_S}\supset\{\sigma\}$ and in this case
$\D_{\sigma_T}\cap\D_{\sigma_S}=\{\sigma\}$ so that the intersection of codimension 1 is irreducible.
\end{proof}

\subsection{}\label{4.6}
Let us compare the results above to the results of Wesbury \cite{W} which were applied by Fung \cite{F} to two-row case. In this subsection we formulate the results which
coincide with our results on two-column case.

For
$T\in\bT_{(n-k,k)^*}$ let $I(T):=I(P_{\sigma_T}):=I(\sigma_T):=\{i\
:\ (i,i+1)\in \sigma_T \}=\{i\ :\ i\in\langle T_1\rangle,\
i+1\in\langle T_2\rangle\}$. For example, take
$$T=\begin{array}{ll}
1&4\\
2&5\\
3&7\\
6&8\\
9&\\
\end{array}$$
 Then $\sigma_{\scriptscriptstyle T}=(3,4)(2,5)(6,7)(1,8)$ and
$I(T)=\{3,6\}.$

Given $i\in \langle T_1\rangle,\ j\in\langle T_2\rangle$ put
$T_{i\leftrightarrows j}$ to be a tableau obtained from $T$ by
interchanging columns of $i$ and $j$ if an array obtained in such a
way is a tableau. Otherwise $T_{i\leftrightarrows j}=\emptyset.$
Note that if $j<i$ then $T_{i\leftrightarrows j}\ne \emptyset$
always. If $j>i$ then $T_{i\leftrightarrows j}\ne\emptyset$ iff
$\sh(\pi_{1,i}(T))=(i-k',k')^*$ where $i-k'> k'+1$ and
$\sh(\pi_{1,j}(T))=(j-k\prpr,k\prpr)^*$ where $j-k\prpr>k\prpr.$

If $i\not\in I(T)$ put
$$u_i(T):=\left\{\begin{array}{ll}
T_{\sigma_T(i)\leftrightarrows i}&{\rm if}\ i,i+1\in \langle
T_2\rangle\\
T_{i+1\leftrightarrows i} &{\rm if}\ i\in \langle T_2\rangle\ {\rm
and}\ i+1\in \langle
T_1\rangle \\
T_{i+1\leftrightarrows \sigma_{T}(i+1)}&{\rm if}\ i,i+1\in \langle
T_1\rangle\ {\rm and}\ \sigma_T(i+1)\ne i+1\\
\emptyset& {\rm if}\ i,i+1\in \langle T_1\rangle\ {\rm and}\
\sigma_T(i+1)=i+1
\end{array}
\right.$$
 Note that in the last case $\sigma_T(i)=i$ as well. And this is
 the only case when $u_i(T)=\emptyset.$

  The picture is much more
 clear on the level of link patterns. Let $T$ be such that
  $i\not\in I(P_T).$ then
 $u_i(P_T)\ne\emptyset$ iff  $\{i,i+1\}\cap\langle P_T\rangle\ne\emptyset.$

 In this case
 \begin{itemize}
 \item{} If $\{i,i+1\}\st \langle\sigma_T\rangle$ then
 $u_i(\sigma_T):=\sigma_{u_i(T)}=(\sigma_T)_{i\leftrightarrows\sigma_T(i+1)}$ and
 $M_{T,u_i(T)}$ consists of $k-2$ loops of length 2 and one loop
 of length 4 with the points
 $\{i,i+1,\sigma_T(i),\sigma_T(i+1)\}.$ There are 3 possible
 situations shown below:
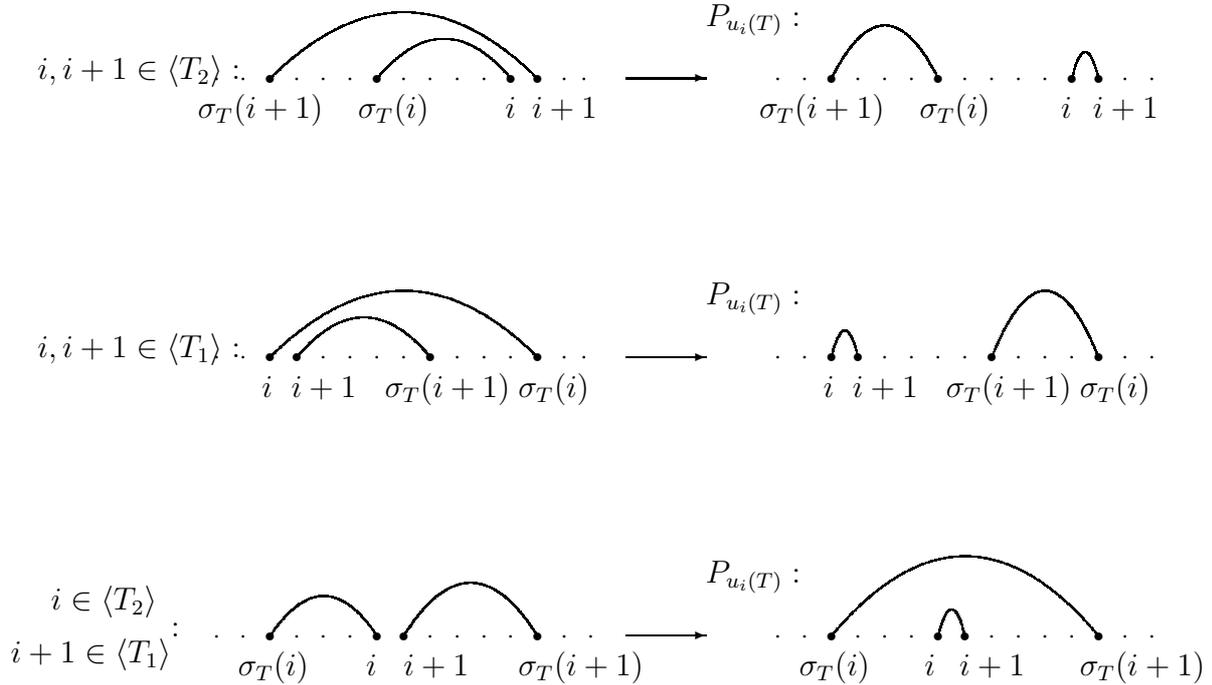
\begin{center}
\begin{picture}(100,100)
\put(-190,60){$i,i+1\in \langle T_2\rangle:$}
\multiput(-123,60)(10,0){2}%
{\circle*{1}}
 \put(-103,60){\circle*{3}}
 \put(-130,45){$\sigma_T(i+1)$}
 \multiput(-93,60)(10,0){3}%
{\circle*{1}}
 \put(-63,60){\circle*{3}}
 \put(-70,45){$\sigma_T(i)$}
\multiput(-53,60)(10,0){4}%
{\circle*{1}}
 \put(-13,60){\circle*{3}}
 \put(-15,45){$i$}
 \put(-3,60){\circle*{3}}
 %\multiput(-47,60)(7.5,0){3}%
%{\circle*{1}}
 \put(-5,45){$i+1$}
\multiput(7,60)(7.5,0){2}%
{\circle*{1}}
 \qbezier(-103,60)(-53,110)(-3,60)
\qbezier(-63,60)(-38,90)(-13,60)

 \put(30,60){\vector(1,0){30}}
 \put(60,80){$P_{u_i(T)}:$}

\multiput(87,60)(10,0){2}%
{\circle*{1}}
 \put(107,60){\circle*{3}}
 \put(80,45){$\sigma_T(i+1)$}
 \put(147,60){\circle*{3}}
 \put(140,45){$\sigma_T(i)$}
\multiput(117,60)(10,0){3}%
{\circle*{1}}
\multiput(157,60)(10,0){4}%
{\circle*{1}}
 \put(197,60){\circle*{3}}
 \put(193,45){$i$}
 \put(207,60){\circle*{3}}
 \put(205,45){$i+1$}
 \multiput(217,60)(10,0){2}%
{\circle*{1}}
 \qbezier(107,60)(127,100)(147,60)
\qbezier(197,60)(202,80)(207,60)

\end{picture}
\end{center}

\begin{center}
\begin{picture}(100,100)
\put(-190,60){$i,i+1\in \langle T_1\rangle:$}
\multiput(-123,60)(10,0){2}%
{\circle*{1}}
 \put(-103,60){\circle*{3}}
 \put(-106,45){$i$}
 \put(-93,60){\circle*{3}}
 \put(-95,45){$i+1$}
 \multiput(-83,60)(10,0){4}%
{\circle*{1}}
 \put(-43,60){\circle*{3}}
 \put(-60,45){$\sigma_T(i+1)$}
\multiput(-33,60)(10,0){3}%
{\circle*{1}}

 \put(-3,60){\circle*{3}}
 %\multiput(-47,60)(7.5,0){3}%
%{\circle*{1}}
 \put(-10,45){$\sigma_T(i)$}
\multiput(7,60)(7.5,0){2}%
{\circle*{1}}
 \qbezier(-103,60)(-53,110)(-3,60)
\qbezier(-93,60)(-68,90)(-43,60)

 \put(30,60){\vector(1,0){30}}
 \put(60,80){$P_{u_i(T)}:$}

\multiput(87,60)(10,0){2}%
{\circle*{1}}
 \put(107,60){\circle*{3}}
 \put(104,45){$i$}
 \put(117,60){\circle*{3}}
 \put(115,45){$i+1$}
\multiput(127,60)(10,0){4}%
{\circle*{1}}
 \put(167,60){\circle*{3}}
 \put(150,45){$\sigma_T(i+1)$}
\multiput(177,60)(10,0){3}%
{\circle*{1}}
 \put(207,60){\circle*{3}}
 \put(200,45){$\sigma_T(i)$}
 \multiput(217,60)(10,0){2}%
{\circle*{1}}
 \qbezier(107,60)(112,80)(117,60)
\qbezier(167,60)(187,110)(207,60)
\end{picture}
\end{center}

\begin{center}
\begin{picture}(100,100)
\put(-185,70){$i\in \langle T_2\rangle$}
 \put(-200,50){$i+1\in \langle T_1\rangle$}
 \put(-140,60){$:$}
\multiput(-123,60)(10,0){2}%
{\circle*{1}}
 \put(-103,60){\circle*{3}}
 \put(-115,45){$\sigma_T(i)$}
\multiput(-93,60)(10,0){3}%
{\circle*{1}}
 \put(-63,60){\circle*{3}}
 \put(-67,45){$i$}
\put(-53,60){\circle*{3}}
 \put(-53,45){$i+1$}
 \multiput(-43,60)(10,0){4}%
{\circle*{1}}
 \put(-3,60){\circle*{3}}
 \put(-10,45){$\sigma_T(i+1)$}
\multiput(7,60)(10,0){2}%
{\circle*{1}}

 \qbezier(-103,60)(-83,90)(-63,60)
\qbezier(-53,60)(-28,100)(-3,60)

 \put(30,60){\vector(1,0){30}}
 \put(60,80){$P_{u_i(T)}:$}

\multiput(87,60)(10,0){2}%
{\circle*{1}}
 \put(107,60){\circle*{3}}
 \put(95,45){$\sigma_T(i)$}
\multiput(117,60)(10,0){3}%
{\circle*{1}}
 \put(147,60){\circle*{3}}
 \put(142,45){$i$}
 \put(157,60){\circle*{3}}
 \put(155,45){$i+1$}
\multiput(167,60)(10,0){4}%
{\circle*{1}}
 \put(207,60){\circle*{3}}
 \put(200,45){$\sigma_T(i+1)$}
 \multiput(217,60)(10,0){2}%
{\circle*{1}}
 \qbezier(107,60)(157,120)(207,60)
\qbezier(147,60)(152,80)(157,60)
\end{picture}
\end{center}
 \item{} If $\{i,i+1\}\not\st\langle\sigma_T\rangle$ then
\begin{itemize}
\item{} If $i\in\langle\sigma_T\rangle$ and $i+1\not\in \langle\sigma_T\rangle$ then $i\in\langle
T_2\rangle,$ $i+1\in\langle T_1\rangle$ and
$u_i(\sigma_T):=\sigma_{u_i(T)}=(\sigma_T)_{\sigma_T(i)\rightarrow
i+1}$, so that $M_{T,u_i(T)}$ consists of $k-1$ loops of length 2
and one interval of length 2 consisting of arcs
$(\sigma_T(i),i),\ (i,i+1)$ as it is shown at the picture below:
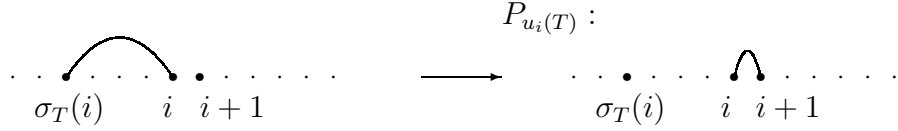
\begin{center}
\begin{picture}(100,100)
%\put(-185,70){$i\in \langle T_2\rangle$}
% \put(-200,50){$i+1\in \langle T_1\rangle$}
 %\put(-140,60){$:$}
\multiput(-123,60)(10,0){2}%
{\circle*{1}}
 \put(-103,60){\circle*{3}}
 \put(-115,45){$\sigma_T(i)$}
\multiput(-93,60)(10,0){3}%
{\circle*{1}}
 \put(-63,60){\circle*{3}}
 \put(-67,45){$i$}
\put(-53,60){\circle*{3}}
 \put(-53,45){$i+1$}
 \multiput(-43,60)(10,0){5}%
{\circle*{1}}
\qbezier(-103,60)(-83,90)(-63,60)

 \put(30,60){\vector(1,0){30}}
 \put(60,80){$P_{u_i(T)}:$}

\multiput(87,60)(10,0){2}%
{\circle*{1}}
 \put(107,60){\circle*{3}}
 \put(95,45){$\sigma_T(i)$}
\multiput(117,60)(10,0){3}%
{\circle*{1}}
 \put(147,60){\circle*{3}}
 \put(142,45){$i$}
 \put(157,60){\circle*{3}}
 \put(155,45){$i+1$}
\multiput(167,60)(10,0){5}%
{\circle*{1}}
 \qbezier(147,60)(152,80)(157,60)
\end{picture}
\end{center}
\item{} If $i+1\in\langle\sigma_T\rangle$ and $i\not\in
\langle\sigma_T\rangle$ then $i,i+1\in\langle T_1\rangle$ and
$u_i(\sigma_T):=\sigma_{u_i(T)}=(\sigma_T)_{\sigma_T(i+1)\rightarrow
i}$ so that $M_{T,u_i(T)}$ consists of $k-1$ loops of length 2 and
an interval of length 2 consisting of arcs $(i,i+1),\
(i+1,\sigma_T(i+1))$ as it is shown at the picture below:
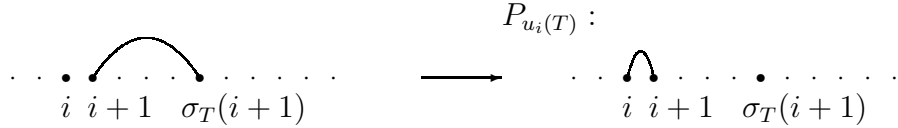
\begin{center}
\begin{picture}(100,100)
%\put(-185,70){$i\in \langle T_2\rangle$}
% \put(-200,50){$i+1\in \langle T_1\rangle$}
 %\put(-140,60){$:$}
\multiput(-123,60)(10,0){2}%
{\circle*{1}}
 \put(-103,60){\circle*{3}}
 \put(-105,45){$i$}
\put(-93,60){\circle*{3}}
 \put(-95,45){$i+1$}
\multiput(-83,60)(10,0){3}%
{\circle*{1}}
 \put(-53,60){\circle*{3}}
 \put(-60,45){$\sigma_T(i+1)$}
 \multiput(-43,60)(10,0){5}%
{\circle*{1}} \qbezier(-93,60)(-73,90)(-53,60)

 \put(30,60){\vector(1,0){30}}
 \put(60,80){$P_{u_i(T)}:$}

\multiput(87,60)(10,0){2}%
{\circle*{1}}
 \put(107,60){\circle*{3}}
 \put(105,45){$i$}
 \put(117,60){\circle*{3}}
 \put(115,45){$i+1$}
\multiput(127,60)(10,0){3}%
{\circle*{1}}
 \put(157,60){\circle*{3}}
 \put(150,45){$\sigma_T(i+1)$}
\multiput(167,60)(10,0){5}%
{\circle*{1}}
 \qbezier(107,60)(112,80)(117,60)
\end{picture}
\end{center}
\end{itemize}
\end{itemize}
Note that in both cases we get by Theorem \ref{4.5} that
$\codim_{\V_T}(\V_T\cap\V_{u_i(T)})=1.$ Moreover, if $i\not\in I(T)$
then $u_i(T)$ is the only tableau in $\{S\in\bT_{(n-k,k)^*}\ :\ i\in
I(S)\}$ such that $\codim_{\V_T}(\V_T\cap\V_{S})=1.$

Further, by \cite[7]{W} one can define inner product on Templerley-Lieb
algebra of $\{P_T\ :\ T\in\bT_{(n-k,k)}\}$ where link pattern
$P_T$ for $T\in\bT_{(n-k,k)}$ is defined exactly in the same way
as it is defined here (that is $P_T:=P_{T^t}.$). Then for two
basis elements $P_T,P_S\ :\ T,S\in\bT_{(n-k,k)}$ one has either
$<P_T,P_S>=\delta^s$ where $0\leq s\leq k-1$ or $<P_T,P_S>=0.$ The
$W$ graph of link patterns with vertices $\{P_T\ :\
T\in\bT_{(n-k,k)}\}$ each labeled by the set $I(T)$ (where $i\in
I(T)$ if $i$ is in the first row of $T$ and $i+1$ in the second
row of $T$, that is $I(T):=I(T^t)$) and edges connecting $P_T$ and
$P_S$ if $<P_T,P_S>=\delta^{k-1}$ is a graph of Kazhdan-Lusztig
type. Moreover, translating the result to our language we get
$<P_{S^t},P_{T^t}>=\delta^{k-1}$ iff $M_{S,T}$ is an even meander
with $k-1$ loops.

Comparing to our results we get that for $S,T\in\bT_{(n-k,k)^*}$
one has $\codim_{\F_S}(\F_T\cap\F_S)=1$ iff
$<P_{S^t},P_{T^t}>=\delta^{k-1}.$ Thus, in particular $W$ graph of
link patterns with vertices $\{P_T\ :\ T\in\bT_{(n-k,k)^*}\}$
labeled by the set $I(T)$ and edges connecting $P_T$ and $P_S$ if
$\codim_{\V_S}(\V_T\cap\V_S)=1$ coincides with the graph of
corresponding Templerley-Lieb algebra and respectively it is of
Kazhdan-Lusztig type.
As it was shown by F. Fung $\codim_{\F_{S^t}}(\F_{T^t}\cap\F_{S^t})=1$ exactly
in the same cases.

\subsection{}\label{4.7}
More generally, in \cite{W} and \cite{G-L} it was shown that for
$P_S,P_T\ :\ S,T\in\bT_{(n-k,k)}$
$$<P_S,P_T>=\left\{\begin{array}{ll}
\delta^r & {\rm if}\ M_{S,T}\ {\rm is\ an\ even\ meander\ with}\ r\
{\rm loops};\\
0 & {\rm if}\ M_{S,T}\ {\rm is\ odd}.\\
\end{array}\right.$$
They showed this theorem using the proposition that for $P_S,P_T\
:\ S,T\in\bT_{(n-k,k)}$ such that $i\in I(P_S)$ and $i\not\in
I(P_T)$ one has
$$<P_S,P_T>=\left\{\begin{array}{ll}
\delta^{-1}<P_S,P_{u_i(T)}> & {\rm if}\ u_i(P_T)\ne\emptyset;\\
0 & {\rm otherwise}.\\
\end{array}\right.$$

F. Fung in \cite{F} considered the intersections of $\F_T\ :\
T\in\bT_{(n-k,k)}.$ By \cite[7]{F} for $T,S\in\bT_{(n-k,k)}$ one
has
$$\codim_{\F_S}(\F_T\cap\F_S)=\left\{\begin{array}{ll}
k-r &{\rm if}\ <P_S,P_T>=\delta^r\\
\emptyset&{\rm if}\ <P_S,P_T>=0\\
\end{array}\right.$$
He showed this  using the same technique, namely showing that for
$i\in I(P_S),\ i\not\in I(P_T)$ one has
$$\codim_{\F_S}(\F_T\cap\F_S)=
\left\{\begin{array}{ll}
 \codim_{\F_S}(\F_S\cap\F_{u_i(T)})+1 & {\rm if}\ u_i(P_T)\ne\emptyset;\\
\emptyset & {\rm otherwise}.\\
\end{array}\right.$$

However the results on the intersections of higher codimensions in two-column case
are very different from those of  Westbury and Fung . First of all note
that since $\F_T\cap\F_S\ne\emptyset$ for $S,T\in\bT_{(n-k,k)^*}$ it
cannot occur that
$\codim_{\F_T}(\F_T\cap\F_S)=\codim_{\F_{T^t}}(\F_{T^t}\cap\F_{S^t})$
for all $S,T\in\bT_{(n-k,k)^*}$.

It is interesting that the deviations of the codimensions are directed in both sides.
On one hand, if $M_{S,T}$ is odd then
$\codim_{\F_T}(\F_T\cap\F_S)<\codim_{F_{T^t}}(\F_{T^t}\cap\F_{S^t})$
( since  $\F_{T^t}\cap\F_{S^t}=\emptyset$). On the other hand,
for even meanders with less than $k-1$ loops it seems that
$\codim_{\F_T}(\F_T\cap\F_S)\geq\codim_{F_{T^t}}(\F_{T^t}\cap\F_{S^t}).$
We do not prove this here, we only show by the example that
the technique of passing from $M_{T,S}$ to $M_{u_i(T),S}$ explained above does not
work in two-column case. Namely we give an example of
$S,T\in\bT_{(n-k,k)^*}$ such that $i\in I(S)$ and $i\not\in I(T)$
and $u_i(T)\ne\emptyset$ and
$\codim_{\F_S}(\F_T\cap\F_S)>\codim_{\F_S}(\F_{u_i(T)}\cap\F_S)+1:$

Let us consider $S,T\in\bT_6^2$ shown below. Note that $2\in I(S)$
and $2\not\in I(T)$ so that
$$S=\begin{array}{ll}
1&3\\
2&4\\
5& \\
6& \\
\end{array},
\quad
 T=\begin{array}{ll}
  1&5\\
2&6\\
3& \\
4& \\
\end{array},
 \quad u_2(T)=\begin{array}{ll}
 1&3\\
 2&5\\
 4& \\
 6&\\
\end{array}$$
The corresponding meanders are:
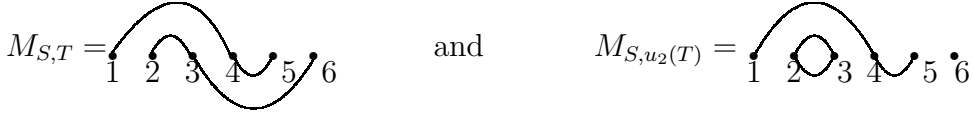
\begin{center}
\begin{picture}(200,80)
\put(-100,40){$M_{S,T}=$}
 \multiput(-60,40)(15,0){6} {\circle*{3}}
 \put(-63,30){1}
 \put(-48,30){2}
 \put(-33,30){3}
 \put(-18,30){4}
 \put(3,30){5}
 \put(18,30){6}

 \qbezier(-61,40)(-33.5,80)(-16,40)
 \qbezier(-46,40)(-38.5,55)(-31,40)
\qbezier(-30,40)(-7.5,0)(16,40)
 \qbezier(-15,40)(-7.5,25)(0,40)
\put(60,40){and}
 \put(120,40){$M_{S,u_2(T)}=$}
 \multiput(180,40)(15,0){6} {\circle*{3}}
\put(177,30){1}
 \put(192,30){2}
 \put(211,30){3}
 \put(222,30){4}
 \put(243,30){5}
 \put(255,30){6}
 \qbezier(180,40)(202.5,80)(225,40)
\qbezier(195,40)(202.5,55)(210,40)
\qbezier(195,40)(202.5,25)(210,40)
\qbezier(240,40)(232.5,25)(225,40)
 \end{picture}
 \end{center}
 By Theorem \ref{4.5}
 $\codim_{\F_S}\F_S\cap\F_{u_2(T)}=1$. On the other hand, a
 straightforward computation shows that $\V_S\cap\V_T$ is
 irreducible and its only component is $\B_{(1,6)(2,5)}$ with
 $\dim\B_{(1,6)(2,5)}=8-4=4$ so that
 $\codim_{\F_S}\F_S\cap\F_T=\codim_{\F_S}\F_S\cap\F_{u_2(T)}+3$.

One can define $I(\sigma)$ and respectively $u_i(\sigma)$ for any $\sigma$ not only
for $\sigma_T$ but in general case $\dim\B_{u_i(\sigma)}\geq \dim \B_\sigma$ and the difference between these two dimensions can be quite big. And this is the reason
of the jumps of codimension of the intersections. But the details of this general
theory are very involved on one hand and  $\F_T\cap\F_S$ of higher codimensions
are not that interesting from the geometric point of view on the other hand, so we stop here.

\subsection{}\label{4.8}
Let us note at the end that there exists a very easy sufficient condition for $\F_T\cap\F_S$ to be reducible. Generalizing this condition we can obtain the necessary
and sufficient condition for the intersection to be irreducible, but again, the formulation of this condition is too involved.

Consider $M_{S,T}$ and put $[i_1,j_1],[i_2,j_2],\ldots$
to be the segments such that $(R_{S,T})_{i_s,j_s}=1$ and $(R_{S,T})_{p,q}=0$
for any $[p,q]\subsetneq [i_s,j_s].$ We order these segments in increasing order
(that is $i_1<i_2<\ldots$) and call them 1-segments of $M_{S,T}.$
\begin{prop} Consider $S,T\in\bT_{(n-k,k)^*}$ Let $[i_1,j_1],\ldots[i_t,j_t]$ be 1-segments of $M_{S,T}.$ If there exists $s\ :\ 1\leq s<t$ such that either $i_{s+1}=j_s$ or $i_{s+1}<j_s$ and $(R_{S,T})_{i_s,j_{s+1}}=1$, then $\F_S\cap\F_T$
is reducible.
\end{prop}
\begin{proof}
First of all note that for a 1-segment $[i_s,j_s]$ there exists a component $\B_\sigma$
of $\ov\B_{\sigma_S}\cap\ov\B_{\sigma_T}$ such that $(i_s,j_s)\in\sigma.$
Indeed, by the definition  of  $[i_s,j_s]$
one has $(i_s,j_s)\prec \sigma_S,\sigma_T.$ Thus, by proposition \ref{3.9a} there exists
a component involution $\sigma$ of $\ov\B_{\sigma_S}\cap\ov\B_{\sigma_T}$ such that $\sigma\succ(i_s,j_s)$. In particular, this provides $(R_\sigma)_{i_s,j_s}\geq 1$.
On the other hand, for any
$[p,q]\subsetneq [i_s,j_s]$ one has
$(R_\sigma)_{p,q}\leq (R_{S,T})_{p,q}=0$ and $(R_\sigma)_{i_s,j_s}\leq (R_{S,T})_{i_s,j_s}=1$. Thus, $(R_\sigma)_{i_s,j_s}= 1$ and $(i_s,j_s)\in\sigma.$

If $i_{s+1}=j_s$ then there exist component involutions $\sigma,\sigma'$ of $\B_T\cap\B_S$ such that
$(i_s,j_s)\in\sigma$ and $(j_s,j_{s+1})\in\sigma'.$ Obviously, $\B_\sigma,\ \B_{\sigma'}$ are two different
components of the intersection.

If $i_{s+1}<j_s$ then there exist component involutions $\sigma,\sigma'$ such that $(i_s,j_s)\in\sigma$ and $(j_s,j_{s+1})\in\sigma'.$ If $\sigma=\sigma'$ we get
$(i_s,j_s),(i_{s+1},j_{s+1})\in\sigma$ so that
$(R_\sigma)_{i_s,j_{s+1}}=2$ which contradicts the condition $(R_{S,T})_{i_s,j_{s+1}}=1$. Thus, again, $\B_\sigma,\ \B_{\sigma'}$ are  two different components of the intersection.
\end{proof}

\bigskip
 \centerline{ INDEX OF NOTATION}
\bigskip

\begin{tabular}{ll}
\ref{1.1}&$\Or_u,\ \gN,\ \bB,\ \B_u,\ \X$\\
\ref{1.2}&$\lambda,\ \lambda^*,\ J(u), \Or_\lambda,\ D_\lambda$\\
\ref{1.3}&$\bS_n^2,\ N_\sigma,\ \B_\sigma,\ P_\sigma$\\
\ref{1.4}&$\ell(P_\sigma),\ c(P_\sigma),\ f_p(P_\sigma),\ f(P_\sigma),\ \bS_n^2(k)$\\
\ref{1.5}&$[i,j],\ R_\sigma,\ \succeq$\\
\ref{1.6}&$\C(\sigma),\ \N(\sigma),\ \D(\sigma)$\\
\ref{1.8}&$R_{\sigma,\sigma'}$\\
\ref{2.1}&$\Bb_{(n-k,k)},\ (i,j),\ ((i,j)),\  q_{(i,j)}(\sigma)$\\
\ref{2.2}&$\V_T,\ \bT_{(n-k,k)^*},\ T_1,\ T_2,\ \langle
T_i\rangle,\ \sigma_T,\ \B_T$\\
\ref{2.3}&$\pi_{i,j},\ R_u,\ R_\sigma,\ \bR_n^2$\\
\ref{3.1}&$ {\bf over}_m(P_\sigma),\ {\bf over}_{(i,j)}(P_\sigma),\ l_{(i,j)}(P_\sigma),\
r_{(i,j)}(P_\sigma),\ c^l_{(i,j)}(P_\sigma),$\\
\ &$\ c^r_{(i,j)}(P_\sigma),\
{\bf under}_{(i,j)}(P_\sigma),\ f'_{(i,j)}(P_\sigma)$\\
\ref{3.2}&$\pi_{i,j}(P_\sigma)$\\
\ref{3.3}&$E(\sigma),\ \sigma_{(i,j)}^- $\\
\ref{3.4}&$\langle P_\sigma\rangle,\ \langle \sigma\rangle,\
\sigma_{i\rightarrow f},\ \sigma_{\curvearrowleft (i,j)},\
           \sigma_{(i,j)\curvearrowright}$\\
\ref{3.5}&$\sigma_{i\leftrightarrows j},\ L_{(i,j)}(\sigma),\ S_{(i,j)\looparrowright}$\\
\ref{3.6}&$Ov_{(i,j)}(\sigma),\
S_{(i,j)\downharpoonleft\upharpoonright}   $\\
\ref{3.8}&$E_{\max}(\sigma)$\\
\end{tabular}

% ----------------------------------------------------------------

%
\end{document}